\documentclass[12pt]{amsart}
\usepackage{amscd,amsmath,amsthm,amssymb,graphics}
\usepackage{amsfonts,amssymb,amscd,amsmath,enumerate,verbatim,calc,latexsym,pstricks,pst-plot,pst-3d,}
\usepackage{amsmath}
\usepackage{amssymb}
\usepackage{amsbsy}
\usepackage{graphicx}
\usepackage{amsfonts}
 \usepackage{multicol}
\usepackage{enumitem}
\usepackage{amscd,amsmath,amsthm,amssymb}
\usepackage{pstricks,pst-plot,pst-3d}
\usepackage{lmodern,pst-node}
\usepackage{pstricks}

\makeatletter

\newcommand{\Rmnum}[1]{\expandafter\@slowromancap\romannumeral #1@}
\makeatother

\usepackage{pgf,tikz}
\usetikzlibrary{arrows}

\def\NZQ{\mathbb}               

\def\RR{{\NZQ R}}

\def\frk{\frak}               

\def\Phi{{\frk n}}
\def\Phi{{\frk N}}
\def\xb{{\mathbf x}}
\def\yb{{\mathbf y}}

\def\wb{{\mathbf w}}
\def\ub{{\mathbf u}}

\def\vb{{\mathbf v}}
\def\opn#1#2{\def#1{\operatorname{#2}}} 
\opn\chara{char} \opn\length{\ell} \opn\pd{pd} \opn\rk{rk}
\opn\projdim{proj\,dim} \opn\injdim{inj\,dim} \opn\rank{rank}
\opn\depth{depth} \opn\grade{grade} \opn\height{height}
\opn\embdim{emb\,dim} \opn\codim{codim}

\opn\Tr{Tr} \opn\bigrank{big\,rank}
\opn\superheight{superheight}\opn\lcm{lcm}
\opn\trdeg{tr\,deg}
\opn\reg{reg} \opn\lreg{lreg} \opn\ini{in} \opn\lpd{lpd}
\opn\size{size}\opn\bigsize{bigsize}
\opn\cosize{cosize}\opn\bigcosize{bigcosize}
\opn\sdepth{sdepth}\opn\sreg{sreg}
\opn\link{link}\opn\fdepth{fdepth}
\opn\index{index}
\opn\index{index}
\opn\indeg{indeg}
\opn\N{N}
\opn\SSC{SSC}
\opn\SC{SC}
\opn\conv{conv}
\opn\div{div} \opn\Div{Div} \opn\cl{cl} \opn\Cl{Cl}
\opn\Spec{Spec} \opn\Supp{Supp} \opn\supp{supp} \opn\Sing{Sing}
\opn\Ass{Ass} \opn\Min{Min}\opn\Mon{Mon} \opn\dstab{dstab} \opn\astab{astab}
\opn\Syz{Syz}
\opn\reg{reg}
\opn\Ann{Ann} \opn\Rad{Rad} \opn\Soc{Soc}
\opn\Im{Im} \opn\Ker{Ker} \opn\Coker{Coker} \opn\Am{Am}
\opn\Hom{Hom} \opn\Tor{Tor} \opn\Ext{Ext} \opn\End{End}
\opn\Aut{Aut} \opn\id{id}

\opn\nat{nat}
\opn\pff{pf}
\opn\Pf{Pf} \opn\GL{GL} \opn\SL{SL} \opn\mod{mod} \opn\ord{ord}
\opn\Gin{Gin} \opn\Hilb{Hilb}\opn\sort{sort}
\opn\initial{init}
\opn\ende{end}
\opn\height{height}
\opn\type{type}
\opn\aff{aff} \opn\con{conv} \opn\relint{relint} \opn\st{st}
\opn\lk{lk} \opn\cn{cn} \opn\core{core} \opn\vol{vol}
\opn\link{link} \opn\star{star}\opn\lex{lex}\opn\Mon{Mon}\opn\Min{Min}
\opn\gr{gr}

\def\pot#1#2{#1[\kern-0.28ex[#2]\kern-0.28ex]}

\opn\dirlim{\underrightarrow{\lim}}
\opn\inivlim{\underleftarrow{\lim}}
\def\Implies{\ifmmode\Longrightarrow \else
        \unskip${}\Longrightarrow{}$\ignorespaces\fi}
\def\implies{\ifmmode\Rightarrow \else
        \unskip${}\Rightarrow{}$\ignorespaces\fi}
\def\iff{\ifmmode\Longleftrightarrow \else
        \unskip${}\Longleftrightarrow{}$\ignorespaces\fi}

\let\:=\colon
\newtheorem{Theorem}{Theorem}[section]
 \newtheorem{Lemma}[Theorem]{Lemma}
 \newtheorem{Corollary}[Theorem]{Corollary}
 \newtheorem{Proposition}[Theorem]{Proposition}
 \newtheorem{Remark}[Theorem]{Remark}
 
 \newtheorem{Example}[Theorem]{Example}
 
 \newtheorem{Definition}[Theorem]{Definition}
 \newtheorem*{Definition*}{Definition}

 \newtheorem*{Conjecture*}{Conjecture}

 \newtheorem{Notation}[Theorem]{Notation}
 \newtheorem{Conventions}[Theorem]{Conventions}

\let\epsilon\varepsilon
\let\kappa=\varkappa
\def\qed{\ifhmode\textqed\fi
      \ifmmode\ifinner\quad\qedsymbol\else\dispqed\fi\fi}
\def\textqed{\unskip\nobreak\penalty50
       \hskip2em\hbox{}\nobreak\hfil\qedsymbol
       \parfillskip=0pt \finalhyphendemerits=0}
\def\dispqed{\rlap{\qquad\qedsymbol}}

\opn\dis{dis}
\def\pnt{{\raise0.5mm\hbox{\large\bf.}}}

\opn\Lex{Lex}

\begin{document}

 \title{Betti numbers of normal edge rings (\bf{I})}

 \author{Zexin Wang}

\address{School  of Mathematical Sciences, Soochow University, 215006 Suzhou, P.R.China}
\email{zexinwang6@outlook.com}

 \author{Dancheng Lu}

\address{School  of Mathematical Sciences, Soochow University, 215006 Suzhou, P.R.China}
\email{ludancheng@suda.edu.cn}

 \begin{abstract}
We introduce  a novel approach named the {\it induced-subgraph approach} for investigating the Betti numbers of  normal edge rings.  Utilizing this approach, we compute all the multigraded Betti numbers of the edge rings associated with two-ear graphs (Definition~\ref{GmDe}) and compact graphs (Definition~\ref{3.1}). In particular, we show that   for two-ear graphs and compact graphs of type 1 or 2, their multigraded Betti numbers are always equal to the top multigraded Betti numbers of some of their induced subgraphs.  In contrast, some of the multigraded Betti numbers of  compact graphs  of type 3 are not the top multigraded Betti numbers of any of their induced subgraphs. We expect that our approach may apply to many other normal edge rings.

 \end{abstract}


 \maketitle

\section{Introduction}

Let $G$ be a simple graph with vertex set
$V(G)=\{x_1,\ldots,x_n\}$ and edge set
$E(G)=\{e_1,\ldots,e_m\}$, and let $\mathbb{K}$ be a field.
The \emph{edge ring} of $G$ is the toric ring
\[
\mathbb{K}[G]
=
\mathbb{K}[x_e\mid e\in E(G)]
\subseteq
\mathbb{K}[x_1,\ldots,x_n],
\]
where
\(x_e:=\prod_{x_i\in e}x_i\)
for every $e\in E(G)$.
Let
\(\mathbb{K}[E(G)]
=
\mathbb{K}[e_1,\ldots,e_m]\)
be the polynomial ring whose variables correspond to the edges of $G$.
The homomorphism
\[
\phi:\mathbb{K}[E(G)]\longrightarrow\mathbb{K}[G],
\qquad
\phi(e_i)=x_{e_i},
\]
is surjective, and its kernel
\[
I_G:=\ker(\phi)
\]
is the \emph{toric ideal} of $G$. Thus
\[
\mathbb{K}[G]\cong\mathbb{K}[E(G)]/I_G.
\]
It is well known that $I_G$ is binomial. Villarreal \cite[Proposition 10.1.10]{V} gives a combinatorial description of its Graver basis, while Tatakis and Thoma \cite{TT} characterize its universal Gr\"obner basis.

We use both the natural multigrading and the standard grading.
Let $\varepsilon_1,\ldots,\varepsilon_n$ be the standard basis of
$\mathbb{Z}^n$, and give $\mathbb{K}[x_1,\ldots,x_n]$ the natural
$\mathbb{Z}_{\geq0}^n$-grading defined by $\deg(x_i)=\varepsilon_i$.
We write $\mathbb{Z}_{\geq0}^{V(G)}$ for the set of all functions from
$V(G)$ to $\mathbb{Z}_{\geq0}$ and identify it with $\mathbb{Z}_{\geq0}^n$ via
$h\mapsto(h(x_1),\ldots,h(x_n))$.
If $e=\{x_i,x_j\}$, set
\[
\deg(e)=\varepsilon_i+\varepsilon_j=\deg(x_e).
\]
Then $\phi$ is multihomogeneous, and hence $I_G$ and
$\mathbb{K}[G]$ inherit a $\mathbb{Z}_{\geq0}^{V(G)}$-grading.
On the other hand, assigning degree $1$ to every edge gives the usual
$\mathbb{Z}_{\geq0}$-grading. According to
\cite[Chapter~IV]{P}, both $I_G$ and $\mathbb{K}[G]$ admit minimal
graded and multigraded free resolutions over $\mathbb{K}[E(G)]$.

Since
\(0\longrightarrow I_G\longrightarrow \mathbb{K}[E(G)]
\longrightarrow \mathbb{K}[G]\longrightarrow 0\)
is a graded and multigraded exact sequence, the Betti numbers of
$I_G$ and $\mathbb{K}[G]$ are related by
$\beta_{i,h}(I_G)=\beta_{i+1,h}(\mathbb{K}[G])$
for all $i\geq0$ and
$h\in\mathbb{Z}_{\geq0}^{V(G)}$, and
$\beta_{i,j}(I_G)=\beta_{i+1,j}(\mathbb{K}[G])$
for all $i,j\geq0$.

If $h=(h_1,\ldots,h_n)\in\mathbb{Z}_{\geq0}^{V(G)}$, set
$|h|=\sum_{k=1}^n h_k$.
Since every edge has standard degree $1$ and multidegree whose
coordinate sum is $2$, a multidegree $h$ corresponds to standard
degree $j$ precisely when $|h|=2j$. Therefore,
\begin{equation}\tag{\dag}\label{dag-betti}
\beta_{i,j}(I_G)
=
\sum_{\substack{
h\in\mathbb{Z}_{\geq0}^{V(G)}\\
|h|=2j}}
\beta_{i,h}(I_G).
\end{equation}

Computing the multigraded, or even the graded, Betti numbers of
$I_G$ is difficult for general graphs, and complete descriptions are
known only for special families. For example, Theorem~5.1 in
\cite{BOV} determines the graded Betti numbers of
$I_{\mathbf{K}_{2,d}}$, where $\mathbf{K}_{2,d}$ is the complete
bipartite graph with bipartition classes
$\{x_1,x_2\}$ and $\{y_1,\ldots,y_d\}$.
Galetto et al.\ \cite{GHK} extended this study to the family
$G_{r,d}$, obtained from $\mathbf{K}_{2,d}$ by adjoining a path of
length $2r-2$ between the two vertices of degree $d$.
More recently, Nandi and Nanduri \cite{NN} investigated the family
$G_{r,s,d}$, obtained by adjoining two such paths.
These are among the few classes of graphs for which the complete
graded Betti numbers have been explicitly described.

In this paper, we develop a unified framework for computing
multigraded Betti numbers of edge rings of graphs satisfying the
\emph{odd cycle condition}. Recall that a simple graph
satisfies the odd cycle condition if any two odd cycles either have a
common vertex or are connected by an edge. For a connected simple
graph $G$, Simis--Vasconcelos--Villarreal \cite{SVV} and
Ohsugi--Hibi \cite{OH} independently showed that this condition is
equivalent to the normality of $\mathbb{K}[G]$.
In particular, the graphs
$\mathbf{K}_{2,d}$, $G_{r,d}$, and $G_{r,s,d}$ belong to this class.

The starting point of our method is the following  property of
multigraded Betti numbers. For
$h\in\mathbb{Z}_{\geq0}^{V(G)}$, let
\[
\operatorname{supp}(h)=\{x\in V(G)\mid h(x)\neq0\},
\qquad
H=G[\operatorname{supp}(h)],
\]
where $H$ is the induced subgraph of $G$ on this support.
By Lemma~\ref{start},
\[
\beta_{i,h}(\mathbb{K}[G])
=
\beta_{i,h}(\mathbb{K}[H])
\qquad\text{for all }i\geq0.
\]
Thus a multigraded Betti number of $\mathbb{K}[G]$ can be computed
inside the induced subgraph determined by its support, without changing
its value.

The key advantage of this reduction is that the induced subgraph
$H=G[\operatorname{supp}(h)]$ may have projective dimension much
smaller than that of $G$. If
\[
\beta_{i,h}(\mathbb{K}[G])
=
\beta_{i,h}(\mathbb{K}[H])\neq0,
\]
then necessarily
\[
i\leq \operatorname{pdim}(\mathbb{K}[H]).
\]
Hence, when $\operatorname{pdim}(\mathbb{K}[H])=i$, the Betti number
$\beta_{i,h}(\mathbb{K}[H])$ is a top multigraded Betti number of
$\mathbb{K}[H]$; when
$\operatorname{pdim}(\mathbb{K}[H])=i+1$, it is a second-top
multigraded Betti number. Thus, for suitable support-induced
subgraphs, the computation of general multigraded Betti numbers can
be reduced to the last one or two steps of their minimal free
resolutions.

This reduction is particularly effective for graphs satisfying the
odd cycle condition. Since this condition is inherited by induced
subgraphs, every induced subgraph $H$ considered above also satisfies
the odd cycle condition. Consequently, $\mathbb{K}[H]$ is normal by
\cite[Corollary~2.3]{OH} and hence, by Hochster's theorem,
Cohen--Macaulay; see, for example, \cite[Theorem~6.10]{BG}. Local duality then relates the top and
second-top multigraded Betti numbers to the minimal generators and
first multigraded Betti numbers of the canonical module; see
Lemma~\ref{canonical}. The canonical module, in turn, has a
combinatorial description in terms of the relative interior of the
edge cone \cite[Theorem~6.31]{BG}, while its first multigraded Betti
numbers can be computed using Lemma~\ref{formula}.

For a fixed homological degree $i$, one of the key ingredients in our
approach is to identify those induced subgraphs whose top multigraded
Betti numbers can contribute to $\beta_i(\mathbb{K}[G])$. This leads
naturally to the notion of a \emph{top-Betti graph}. More precisely, we
call a simple graph $H$ a \emph{top-Betti graph} if every nonzero top
multigraded Betti number has full support, that is,
\[
\beta_{p,h}(\mathbb{K}[H])\neq 0
\quad\Longrightarrow\quad
\operatorname{supp}(h)=V(H)
\]
for every $h\in\mathbb{Z}_{\geq 0}^{V(H)}$, where
$p=\operatorname{pdim}(\mathbb{K}[H])$.

Since a nonzero top multigraded Betti number contributed by a
top-Betti induced subgraph $H$ has support exactly $V(H)$, 
distinct top-Betti induced subgraphs contribute in disjoint
multidegrees.  Therefore the sum of these contributions gives a lower bound
for $\beta_i(\mathbb{K}[G])$.

Another key ingredient is an upper bound obtained from a suitable
initial ideal:
\[
\beta_i(\mathbb{K}[G])
\leq
\beta_i\bigl(
\mathbb{K}[E(G)]/\operatorname{in}_{<}(I_G)
\bigr).
\]
We compare the  above lower bound with this
upper bound. If they agree, then the top-Betti contributions
determine all multigraded Betti numbers in homological degree $i$.
If a gap remains, further contributions must be taken into account.
In this case, we consider induced subgraphs $H$ satisfying
\[
\operatorname{pdim}(\mathbb{K}[H])=i+1,
\]
whose second-top multigraded Betti numbers contribute in homological
degree $i$. This leads naturally to the notion of a
\emph{second-Betti graph}, which is introduced formally in Section~3.

We call this method the \emph{induced-subgraph approach}.
Multi-path graphs, including $\mathbf{K}_{2,d}$, $G_{r,d}$, and
$G_{r,s,d}$, are studied using this approach in \cite{WL1}. In the
present paper, we concentrate on two-ear graphs and compact graphs.

The first family consists of the graphs $\mathbf{G}_m$, introduced in
\cite{HHKO}; see Definition~\ref{GmDe}. Owing to their shape, we call
them \emph{two-ear graphs}. We show that every nonzero multigraded
Betti number of the edge ring of a two-ear graph is obtained from a
top multigraded Betti number of a suitable induced subgraph.
Moreover, the relevant top-Betti induced subgraphs belong to only
three explicit families: two-ear graphs, one-ear graphs, and complete
bipartite graphs of the form $\mathbf{K}_{2,d}$. As a consequence,
all multigraded Betti numbers of the edge ring of a two-ear graph are
either $0$ or $1$. In particular, the top-Betti contributions attain
the upper bounds coming from the initial ideal, so no second-top
contributions are needed for this family.

The second family consists of compact graphs. Following \cite{WL}, a
\emph{compact graph} is a connected simple graph of minimum degree at
least $2$ which contains no even cycle and satisfies the odd cycle
condition. A vertex of a compact graph is called \emph{big} if its
degree is greater than $2$. By \cite[Lemma~3.10]{WL}, a compact graph
has at most three big vertices. We say that a compact graph is of
type $r$ if it has exactly $r$ big vertices, where $0\leq r\leq3$.
The case $r=0$ is an odd cycle and is trivial for our purposes. 
For compact graphs of type $1$ or $2$, as for two-ear graphs, the
top-Betti contributions attain the upper bounds provided by a suitable
initial ideal, so no second-top contributions are needed. Every
nonzero multigraded Betti number is obtained from a top multigraded
Betti number of a compact induced subgraph, and all multigraded Betti
numbers are either $0$ or $1$. For compact graphs of type $3$,
however, the top-Betti contributions do not in general attain these
upper bounds. To fill the gap, certain second-top multigraded Betti
numbers of type-$3$ compact induced subgraphs are also required, and
multigraded Betti numbers equal to $2$ may occur. Thus second-top
contributions are genuinely necessary in this case.

This distinction suggests the following natural question: \emph{which graphs
satisfying the odd cycle condition have the property that all
multigraded Betti numbers of their edge rings can be recovered from
top multigraded Betti numbers of induced subgraphs?}  The results of this
paper, together with \cite{WL1}, show that two-ear graphs, compact
graphs of type $1$ or $2$, and multi-path graphs have this property,
whereas compact graphs of type $3$ show that it does not hold in
general.

The paper is organized as follows. In Section~2 we recall basic
definitions and preliminary results. In Section~3 we introduce and study top-Betti and second-Betti graphs. Section~4 is
devoted to the multigraded Betti numbers of two-ear graphs, together
with their one-ear and $\mathbf{K}_{2,d}$ induced subgraphs.
In Section~5 we study the canonical modules and top supports of compact
graphs. Finally, in Section~6 we apply the induced-subgraph approach
to determine the multigraded Betti numbers of compact graphs.

\section{Preliminaries}

In this section, we briefly review the notation and fundamental facts that will be utilized later on.

\subsection{Betti numbers}
Let $\mathbb{K}$ be a field, and let $R:=\mathbb{K}[X_1,\ldots,X_n]$ be the polynomial ring in variables $X_1,\ldots,X_n$ that are $\mathbb{Z}^m$-graded. Here, we do not require that $m=n$.   For a finitely generated $\mathbb{Z}^m$-graded $R$-module $M$, there exists a minimal multigraded free resolution of $M$ that has the form:
\begin{equation*}0\rightarrow\underset{h\in\mathbb{Z}^m}{\bigoplus}R[-h]^{\beta_{p,h}(M)}\rightarrow \cdots \rightarrow \underset{h\in\mathbb{Z}^m}{\bigoplus}R[-h]^{\beta_{1,h}(M)}\rightarrow\underset{h\in\mathbb{Z}^m}{\bigoplus}R[-h]^{\beta_{0,h}(M)}\rightarrow M \rightarrow 0.\end{equation*}
Note that $R[-h]$ is the cyclic free $R$-module generated in multidegree $h$.
The {\it projective dimension} of $M$ is defined to be  $\mathrm{pdim}\,(M):={\mbox{max}}\,\{i\mid \beta_{i,\,h}(M)\neq 0 \mbox{ for some } h \}.$
The number $\beta_{i,h}(M):={\rm{dim}}_{\mathbb{K}}\mathrm{Tor}_i^R(M,\mathbb{K})_h$ is called the $(i,h)$-th {\it multigraded Betti number} of $M$ and $\beta_{i}(M):=\sum_{h\in \mathbb{Z}^m}\beta_{i,h}(M)$ is called the $i$-th {\it total Betti number} of $M$. When $R$ and $M$ are $\mathbb{Z}$-graded, the number
$\beta_{i,j}(M):=\dim_{\mathbb{K}}\mathrm{Tor}_i^R(M,\mathbb{K})_j$
is called the $(i,j)$-th {\it graded Betti number} of $M$.
The {\it regularity} of $M$ is defined by \(\mathrm{reg}(M):=\max\{j-i\mid \beta_{i,j}(M)\neq 0\}.\)

Denote by $p$ the projective dimension of $M$. Subsequently, $\beta_p(M)$ is referred to as the {\it top total Betti number} of $M$, while $\beta_{p,h}(M)$ for $h\in \mathbb{Z}^m$ are called the {\it top multigraded Betti numbers} of $M$. Similarly, $\beta_{p-1}(M)$ is labeled as the {\it second-top total Betti number} of $M$, and $\beta_{p-1,h}(M)$ for $h\in \mathbb{Z}^m$ are called the {\it second-top multigraded Betti numbers} of $M$.

To obtain upper bounds for the Betti numbers of edge rings, we will use the following
well-known result (see, for instance, \cite[Corollary 3.3.3]{HH}).

\begin{Lemma}\label{total}
Let $I$ be a graded ideal of the polynomial ring  $R=\mathbb{K}[X_1, \ldots, X_n]$, and let $<$  be any monomial order on $R$. Then for all $i, j \geq 0$, we have
 $$\beta_{i,j}(R/I) \leq \beta_{i,j}(R/\mathrm{in}_<(I)).$$
 In particular,  $$\beta_{i}(R/I) \leq \beta_{i}(R/\mathrm{in}_<(I)) .$$
\end{Lemma}

\subsection{Simple graphs and minimal cycles}
By definition, a simple graph is a finite graph devoid of loops and multiple edges.
Let $G$ be a simple graph and denote  the  vertex set and edge set of $G$  by $V(G)$ and  $E(G)$ respectively.   Recall that a \emph{walk} of $G$ of length $q$ is a sequence $i_0,i_1,\ldots,i_q$ of vertices such that $\{i_{j-1},i_j\}\in E(G)$ for $1\leq j\leq q$.
  The graph $G$ is \emph{connected} if, for any two vertices $i_j$ and $i_k$ of $G$, there
 is a walk between $i_j$ and $i_k$.

 A \emph{cycle} of $G$ of length $q$ is a subgraph $C$ of $G$ such that
 $$E(C)=\{\{i_1,i_2\},\{i_2,i_3\},\ldots, \{i_{q-1},i_q\},\{i_q,i_1\}\}$$
 where $i_1,i_2,\ldots,i_q$ are vertices of $G$ and where $i_j\neq i_k$ if $j\neq k$.
 A cycle is called \emph{even} (resp. \emph{odd}) if $q$ is even (resp. odd). For a subset $W$ of $V(G)$, the {\it induced subgraph} $G[W]$ is the graph with vertex set $W$ and for every pair $j,k\in W$, they are adjacent in $G[W]$ if and only if they are adjacent in $G$. For a vertex $j\in V(G)$, let $G \setminus j$ be the induced subgraph $G[V(G) \setminus \{j\}]$. Throughout this paper, all induced subgraphs under consideration are
 assumed to be nonempty unless otherwise stated. A cycle of $G$ is called an \emph{induced cycle} if it is an induced subgraph of $G$. We also call an induced cycle of $G$ to be a \emph{minimal cycle} of $G$, and use $\mathfrak{t}(G)$ to denote the number of minimal cycles of $G$. A non-empty subset $T \subseteq V(G)$ is called an \textit{independent set} of $G$ if $\{j,k\} \not\in E(G)$ for any $j,k \in T$.  For an independent set $T \subseteq V(G)$, we use $N_G(T)$  (or simply $N(T)$ for brevity) to  denote the neighborhood set of $T$, namely, $$N_G(T):=\{j\in V(G)\mid \{j,k\}\in E(G) \mbox{ for some  } k\in T\}.$$ For a vertex $x\in V(G)$, its degree is defined by $\deg_G(x):=|N_G(\{x\})|$; when the ambient graph is clear, we write simply $\deg(x)$.

The following two lemmas are important, especially the former, which serves as the cornerstone
of the induced-subgraph approach. The first lemma is a special case of \cite[Corollary~2.5]{OHH},
and the second one is also well known in the literature; we record them here, and for the second
we include a proof for completeness.

\begin{Lemma}\label{start}
Let $H$ be an induced subgraph of a simple graph $G$. Then, for every
$h\in \mathbb{Z}_{\geq 0}^{V(G)}$ with $\mathrm{supp}(h)\subseteq V(H)$, we have
\[
\beta_{i,h}(\mathbb{K}[H])=\beta_{i,h}(\mathbb{K}[G])
\qquad \text{for all } i\geq 0.
\]
\end{Lemma}

Let $c_0(G)$ denote the number of bipartite connected components of a simple graph $G$.

\begin{Lemma} \label{pd}
Let $G$ be a simple graph with vertex set $V(G)$ and edge set $E(G)$. If $G$ satisfies the
odd cycle condition, then
\[
\mathrm{pdim}(\mathbb{K}[G])=|E(G)| - |V(G)| + c_0(G).
\]
\end{Lemma}

\begin{proof}
First consider the case when $G$ is connected. By \cite[Proposition 13.15]{S} and
\cite[Theorem 6.3.5]{BH}, $\mathbb{K}[G]$ is Cohen-Macaulay.
It follows from \cite[Corollary 10.1.21]{V} that
\[
\mathrm{depth}(\mathbb{K}[G]) = \mathrm{dim}(\mathbb{K}[G]) = \left\{
\begin{array}{ll}
|V(G)| - 1, & \text{if $G$ is bipartite;} \\
|V(G)|,     & \text{otherwise.}
\end{array}
\right.
\]
According to the Auslander-Buchsbaum formula, as stated, for instance, in
\cite[Corollary A.4.3]{HH}, we have
\[
 \mathrm{pdim}(\mathbb{K}[G]) = \left\{
\begin{array}{ll}
|E(G)| - |V(G)| + 1, & \text{if $G$ is bipartite;} \\
|E(G)| - |V(G)|,     & \text{otherwise.}
\end{array}
\right.
\]
Next, assume $G$ has connected components $H_1,\ldots,H_k$. Then, by \cite[Lemma 3.8]{HBO}, we have
\[
\mathrm{pdim}(\mathbb{K}[G]) = \mathrm{pdim}(\mathbb{K}[H_1]) + \cdots + \mathrm{pdim}(\mathbb{K}[H_k]).
\]
Now the result follows by combining the above two equalities.
\end{proof}

\section{Top-Betti graphs and second-Betti graphs}

In this section, we study top-Betti and second-Betti graphs, which
play a central role in our approach. We first record the definition
of a top-Betti graph.

\begin{Definition} \label{BettiD}\em
Let \( G \) be a simple graph, and let \( p = \mathrm{pdim}(\mathbb{K}[G]) \) denote the projective dimension of \( \mathbb{K}[G] \). We call \( G \) a \emph{top-Betti graph} if whenever \( \beta_{p,h}(\mathbb{K}[G]) \neq 0 \) for some \( h \in \mathbb{Z}_{\geq 0}^{V(G)} \), it follows that \( \operatorname{supp}(h) = V(G) \). Recall that $\operatorname{supp}(h)=\{x\in V(G)\mid h(x)\neq 0\}$.

\end{Definition}

Combining Definition~\ref{BettiD} with Lemma~\ref{start}, we obtain the following characterization.

\begin{Lemma}\label{lem:top-Betti-pdim}
Let $G$ be a simple graph and set $p=\mathrm{pdim}(\mathbb{K}[G])$. Then $G$ is a top-Betti graph if and only if for every vertex $x\in V(G)$ one has
\[
  \mathrm{pdim}(\mathbb{K}[G\setminus x])<\mathrm{pdim}(\mathbb{K}[G])=p.
\]
\end{Lemma}

By \cite[Lemma 3.8]{HBO} and Lemma~\ref{lem:top-Betti-pdim}, we obtain the following corollary.
\begin{Corollary}\label{top-components} A simple graph $G$ is a top-Betti graph if and only if all its  connected components are top-Betti graphs.
\end{Corollary}

We next present a combinatorial characterization of a top-Betti graph for connected graphs satisfying the odd cycle condition.
Recall that $c_0(H)$ denotes the number of bipartite connected components of a graph $H$.

\begin{Proposition}\label{Betti}
Let $G$ be a connected simple graph satisfying the odd cycle condition. Then $G$ is a top-Betti graph if and only if for every $x\in V(G)$,
$c_0(G\setminus x) < \deg (x)$ if $G$ is bipartite, and $c_0(G\setminus x) < \deg(x) - 1$ if $G$ is non-bipartite.
\end{Proposition}

\begin{proof}
Set $p=\mathrm{pdim}(\mathbb{K}[G])$. By Lemma~\ref{pd}, we get
\[
  p = |E(G)|-|V(G)|+c_0(G),
\]
and
\[
  \mathrm{pdim}(\mathbb{K}[G\setminus x])
  = |E(G\setminus x)|-|V(G\setminus x)|+c_0(G\setminus x).
\]
Since $|E(G\setminus x)|=|E(G)|-\deg (x)$ and $|V(G\setminus x)|=|V(G)|-1$, it follows that
\[
  \mathrm{pdim}(\mathbb{K}[G\setminus x])
  = |E(G)|-|V(G)|+1-\deg (x) + c_0(G\setminus x).
\]
Therefore, for each $x\in V(G)$ we have
\[
  \mathrm{pdim}(\mathbb{K}[G\setminus x])<p
  \;\Longleftrightarrow\;
  |E(G)|-|V(G)|+1-\deg (x) + c_0(G\setminus x)
  < |E(G)|-|V(G)|+c_0(G),
\]
which is equivalent to
\[
  c_0(G\setminus x) < \deg (x) + c_0(G) - 1.
\]
Since $G$ is connected, we have
\[
  c_0(G)=
  \begin{cases}
    1, & \text{if $G$ is bipartite},\\
    0, & \text{if $G$ is non-bipartite},
  \end{cases}
\]
which yields the stated piecewise bounds.
\end{proof}

As an immediate consequence of Proposition~\ref{Betti}, we obtain the following simple degree constraint on top-Betti graphs, which allows us to determine that many induced subgraphs under consideration are not top-Betti graphs.

\begin{Corollary}\label{top-Betti at least two}
Let $G$ be a simple graph that satisfies the odd cycle condition.
If $G$ is a top-Betti graph, then every vertex of $G$ has degree at least two.
\end{Corollary}

\begin{proof}
	If $G$ had an isolated vertex $x$, then
	$\mathbb{K}[G\setminus x]=\mathbb{K}[G]$, contradicting
	Lemma~\ref{lem:top-Betti-pdim}. Thus $G$ has no isolated vertices.
	
	Let $G_1,\ldots,G_s$ be the connected components of $G$. By
	Corollary~\ref{top-components}, each $G_\ell$ is a top-Betti graph.
	Moreover, each $G_\ell$ satisfies the odd cycle condition, so
	Proposition~\ref{Betti} applies to every $G_\ell$.
	
	Fix $1\leq \ell\leq s$ and suppose that
	$\deg_{G_\ell}(x)=1$ for some $x\in V(G_\ell)$. If $G_\ell$ is
	bipartite, then $G_\ell\setminus x$ is connected and bipartite, so
	\[
	c_0(G_\ell\setminus x)=1=\deg_{G_\ell}(x),
	\]
	contradicting Proposition~\ref{Betti}. If $G_\ell$ is non-bipartite,
	then $G_\ell\setminus x$ is connected and non-bipartite, so
	\[
	c_0(G_\ell\setminus x)=0=\deg_{G_\ell}(x)-1,
	\]
	again contradicting Proposition~\ref{Betti}. Hence every vertex of
	$G$ has degree at least $2$.
\end{proof}

We use $\sqcup$ to denote disjoint union. We now give two examples illustrating the above notions and results.

\begin{Example}\em
Let $\mathbf{K}_{m,n}$ denote the complete bipartite graph with bipartition
$U\sqcup W$, where $|U|=m$ and $|W|=n$, in which every vertex of $U$
is adjacent to every vertex of $W$. If either $m=1$ or $n=1$, then  $\mathbf{K}_{m,n}$ is not a top-Betti graph by Corollary~\ref{top-Betti at least two}. If both $m\geq 2$ and $n\geq 2$, then $c_0(\mathbf{K}_{m,n}\setminus x)=1<\deg (x)$ for any $x\in V(\mathbf{K}_{m,n})$, and hence $\mathbf{K}_{m,n}$ is a top-Betti graph by Proposition~\ref{Betti}.
\end{Example}

A graph $G$ with at least three vertices is called \emph{$2$-connected} if
	$G\setminus x$ is connected for every $x\in V(G)$.

\begin{Example}\em 	
	Let $G$ be a $2$-connected graph satisfying the odd cycle condition.
	Then $G$ is a top-Betti graph if and only if either $G$ is bipartite or
	$G\setminus x$ is non-bipartite for every vertex $x\in V(G)$ of degree $2$.
	
Indeed, $2$-connectedness implies that $G\setminus x$ is connected and
$\deg_G(x)\geq 2$ for every $x\in V(G)$. If $G$ is bipartite, then
\[
c_0(G\setminus x)=1<\deg_G(x)
\]
for every $x\in V(G)$, so Proposition~\ref{Betti} shows that $G$ is top-Betti.
Assume that $G$ is non-bipartite. Since $G\setminus x$ is connected,
$c_0(G\setminus x)$ is equal to $1$ precisely when $G\setminus x$ is bipartite,
and is equal to $0$ otherwise. Hence the inequality
\[
c_0(G\setminus x)<\deg_G(x)-1
\]
fails if and only if $c_0(G\setminus x)=1$ and $\deg_G(x)=2$.
Therefore, by Proposition~\ref{Betti}, $G$ is top-Betti if and only if
$G\setminus x$ is non-bipartite for every vertex $x\in V(G)$ of degree $2$. For example, the graph consisting of two triangles sharing a common edge
is a top-Betti graph, since deleting either of its vertices of degree $2$
leaves a triangle.
\end{Example}

In the induced-subgraph approach, the top-Betti contributions in a
fixed homological degree $i$ may not account for all multigraded Betti
numbers in that degree. In this case, we need to consider second-top
multigraded Betti numbers of induced subgraphs of projective dimension
$i+1$. This motivates the following definition.
\begin{Definition}\label{second-BettiD}{\em
Let $G$ be a simple graph with
$p=\operatorname{pdim}(\mathbb{K}[G])\geq 1$.
We call $G$ a \emph{second-Betti graph} if whenever
$\beta_{p-1,h}(\mathbb{K}[G])\neq 0$ for some multidegree $h$, then either
$\mathrm{supp}(h)=V(G)$ or the projective dimension of
$\mathbb{K}[G[\mathrm{supp}(h)]]$ is equal to $p-1$.
A second-Betti graph is termed \emph{proper} if there is a multidegree $h$ such that
$\mathrm{supp}(h)=V(G)$ and $\beta_{p-1,h}(\mathbb{K}[G])\neq 0$.}
\end{Definition}

The following proposition will be used to identify proper second-Betti induced subgraphs of the graphs under consideration.

\begin{Proposition} \label{second-Betti}For any simple graph $G$, if $G$ is a top-Betti graph, then $G$ is a second-Betti graph.
\end{Proposition}
\begin{proof}
	Set $p=\operatorname{pdim}(\mathbb{K}[G])$.
	Let $h$ be a multidegree such that
	$\beta_{p-1,h}(\mathbb{K}[G])\neq 0$. We are done if $\mathrm{supp}(h)=V(G)$. Suppose that $\mathrm{supp}(h)\subsetneqq V(G)$. Let $H$ be the subgraph of $G$ induced on $\mathrm{supp}(h)$. Note that, since $H$ is a proper induced subgraph of $G$, it follows from
Lemma~\ref{lem:top-Betti-pdim} that
\[
\mathrm{pdim}(\mathbb{K}[H])<\mathrm{pdim}(\mathbb{K}[G]),
\]
and thus $\mathrm{pdim}(\mathbb{K}[H])\leq p-1$.
 However, we have $\beta_{p-1, h}(\mathbb{K}[H])=\beta_{p-1, h}(\mathbb{K}[G])\neq 0$. Hence,  $\mathrm{pdim}(\mathbb{K}[H])=p-1$, as desired.
\end{proof}

At the end of this section, we introduce the notion of the top support
of a simple graph, which will play an important role in the subsequent
sections. It records precisely the multidegrees in which the top
multigraded Betti numbers are nonzero.

\begin{Definition}\label{topsupport}
\emph{Let $H$ be a simple graph, and let
$p=\operatorname{pdim}(\mathbb{K}[H])$.
We define
\[
\mathcal{N}_H
:=
\left\{
X^h \;\middle|\;
\beta_{p,h}(\mathbb{K}[H])\neq 0
\right\},
\]
and call $\mathcal{N}_H$ the \emph{top support} of $H$.}
\end{Definition}
In the situations considered below, every nonzero top multigraded
Betti number is equal to $1$. Hence, the top support $\mathcal{N}_H$
completely determines the top multigraded Betti numbers of
$\mathbb{K}[H]$.

\section{Betti numbers of $\mathbf{G}_m$}

In \cite{HHKO}, a specific class of simple graphs, designated as
$\mathbf{G}_m$, was introduced for all $m\geq 1$. It was proved there
that the toric ideal of $\mathbb{K}[\mathbf{G}_m]$ admits a squarefree
initial ideal \cite[Corollary~4.3]{HHKO} whose quotient ring is
Cohen--Macaulay \cite[Lemma~4.1]{HHKO}. We refer to such a graph as a
\emph{two-ear} graph due to its shape.

In this section, we employ the induced-subgraph approach to demonstrate
that the toric ideal of $\mathbf{G}_m$ possesses exactly the same Betti
numbers as this squarefree initial ideal. Furthermore, we present all
of its multigraded Betti numbers.

\begin{Definition}\label{GmDe}\em
Let $m \geq 1$ be an integer. By relabeling the vertices and edges of the graph introduced in \cite[Figure~1]{HHKO}, we define a simple graph $\mathbf{G}_{m}$ with vertex set
\[
V(\mathbf{G}_{m})
= \{v_1,v_2,x_1,x_2,y_1,y_2\} \sqcup \{u_{i} \mid 1 \leq i \leq m-1\}
\]
and edge set
\begin{align*}
E(\mathbf{G}_{m}) = {} &
\big\{\{v_{1},x_{1}\},\{v_{1},x_{2}\},\{x_{1},x_{2}\},\{v_{2},y_{1}\},\{v_{2},y_{2}\},\{y_{1},y_{2}\},\{x_{1},y_{1}\}\big\} \\
& \sqcup \big\{\{v_{1},u_{i}\} \mid 1 \leq i \leq m-1\big\}
   \sqcup \big\{\{v_{2},u_{i}\} \mid 1 \leq i \leq m-1\big\}.
\end{align*}
Let $\mathbf{P}_m$ denote the induced subgraph of $\mathbf{G}_{m}$ obtained by deleting the vertex $x_2$, that is,
\[
\mathbf{P}_m = \mathbf{G}_{m} \setminus x_2.
\]
We call a graph $G$ a {\it two-ear graph of type $m$} if it is isomorphic to $\mathbf{G}_m$.
Similarly, we call a graph $G$ a {\it one-ear graph of type $m$} if it is isomorphic to $\mathbf{P}_m$.
\end{Definition}

\begin{figure}[ht]
\centering
\begin{tikzpicture}

\draw[black, thin] (3,1) -- (3,5) -- (6,3) -- (9,5) -- (9,1)-- cycle;
\draw[black, thin] (6,6) -- (3,5) --(6,7.5) -- (9,5)-- cycle;
\draw[black, thin] (3,1) -- (3,5) -- (1.3,3) --  cycle;
\draw[black, thin] (10.7,3) -- (9,5) -- (9,1) --  cycle;

\filldraw [black] (3,5) circle (2pt);
\filldraw [black] (9,5) circle (2pt);
\filldraw [black] (3,1) circle (1pt);
\filldraw [black] (6,3) circle (1pt);
\filldraw [black] (9,1) circle (1pt);
\filldraw [black] (6,6) circle (1pt);
\filldraw [black] (6,7.5) circle (1pt);
\filldraw [black] (1.3,3)  circle (1pt);
\filldraw [black] (10.7,3)  circle (1pt);

\filldraw [black] (6,4.2) circle (0.5pt);
\filldraw [black] (6,4.5)  circle (0.5pt);
\filldraw [black] (6,4.8)  circle (0.5pt);

\draw (3,5) node[anchor=south]{$v_1$};
\draw (9,5) node[anchor=south]{$v_2$};
\draw (3,1) node[anchor=north]{$x_1$};
\draw (6,3) node[anchor=north]{$u_{m-1}$};
\draw (9,1) node[anchor=north]{$y_1$};
\draw (6,6.5) node[anchor=north]{$u_2$};
\draw (6,7.5) node[anchor=south]{$u_1$};
\draw (1.3,3) node[anchor=east]{$x_2$};
\draw (10.7,3) node[anchor=west]{$y_2$};

\draw (4.2,4) node[anchor=north]{$e_{m-1,1}$};
\draw (7.9,4) node[anchor=north]{$e_{m-1,2}$};
\draw (4.2,5.2) node[anchor=north]{$e_{2,1}$};
\draw (7.8,5.2) node[anchor=north]{$e_{2,2}$};
\draw (4,6) node[anchor=south]{$e_{1,1}$};
\draw (8,6) node[anchor=south]{$e_{1,2}$};
\draw (3,3) node[anchor=west]{$f_1$};
\draw (9,3) node[anchor=east]{$g_1$};
\draw (2.1,1.7) node[anchor=east]{$f_2$};
\draw (9.9,1.7) node[anchor=west]{$g_2$};
\draw (2.1,4.3) node[anchor=east]{$f_3$};
\draw (9.9,4.2) node[anchor=west]{$g_3$};
\draw (6,0.9) node[anchor=north]{$w$};
\end{tikzpicture}
\caption{The two-ear graph $\mathbf{G}_{m}$}
\label{fig3}

\end{figure}

 Let us denote $\mathbf{G}_m$ by $G$. To write out the minimal generators of $\mathrm{in}_{<}(I_G)$,  we label the edges of $G$ as follows: for $1\leq i\leq m-1$, $e_{i,1}=\{v_1,u_{i}\},e_{i,2}=\{v_{2},u_{i}\}$, and $f_1=\{v_1,x_1\}, f_2=\{x_1,x_2\}, f_3=\{v_1,x_2\}, g_1=\{v_2,y_1\}, g_2=\{y_1,y_2\}, g_3=\{v_2,y_2\}, w=\{x_1,y_1\}$.   In \cite[Corollary 4.3 ]{HHKO}, the minimal set of monomial generators of $\mathrm{in}_{<}(I_{G})$  is specified for a particular monomial order $<$.  Note that the monomial order $<$ in this section refers to this order. We record this result in the following lemma for later use.

\begin{Lemma} There is a monomial order $<$ such that the minimal set of monomial generators of $\mathrm{in}_{<}(I_{G})$ is given  by $\mathcal{M}=\mathcal{M}_{1}\sqcup\mathcal{M}_{2}\sqcup\mathcal{M}_{3}\sqcup\mathcal{M}_{4}$, where
\begin{enumerate}
\item[$(i)$] $\mathcal M_1= \left\{e_{i,1}e_{j,2}\ |\   1 \leq i < j \leq m-1 \right\}$,
\item[$(ii)$] $\mathcal M_2= \left\{f_1f_2g_1g_2 \right\}$,
\item[$(iii)$] $\mathcal M_3=\left\{g_{1}f_{2}e_{i,1}\ |\   1 \leq i \leq m-1  \right\}$, and
\item[$(iv)$] $\mathcal M_4=\left\{g_{2}f_{1}e_{i,2}\ |\   1 \leq i \leq m-1 \right\}$.
\end{enumerate}

\end{Lemma}

We order the minimal generators of $\mathrm{in}_{<}(I_G)$ as follows:
first the generators $e_{i,1}e_{j,2}$, $1\leq i<j\leq m-1$, in increasing
order of $i$ and, for each fixed $i$, in increasing order of $j$; then
$g_1f_2e_{1,1}\prec g_1f_2e_{2,1}\prec\cdots\prec g_1f_2e_{m-1,1}$;
followed by $f_1f_2g_1g_2$; and finally
$g_2f_1e_{m-1,2}\prec g_2f_1e_{m-2,2}\prec\cdots\prec g_2f_1e_{1,2}$.
For brevity, we label the terms in this sequence consecutively as
$\alpha_1\prec\cdots\prec\alpha_r$, where
$r=\frac{m(m+1)}{2}$.
\begin{Lemma}\label{Gmin}
	Let $I=\mathrm{in}_{<}(I_G)$ and keep the above order
	$\alpha_1\prec\cdots\prec\alpha_r$ on its minimal generators.
	For $j=2,\ldots,r$, the colon ideal
	\[
	(\alpha_1,\ldots,\alpha_{j-1}):\alpha_j
	\]
	is generated by a regular sequence of length $k_j$, where
	\[
	k_j=
	\begin{cases}
		b-2, & \text{if }\alpha_j=e_{a,1}e_{b,2},\\
		m-2, & \text{if }\alpha_j=g_1f_2e_{a,1},\\
		m-1, & \text{if }\alpha_j=f_1f_2g_1g_2,\\
		m-1, & \text{if }\alpha_j=g_2f_1e_{a,2}.
	\end{cases}
	\]
\end{Lemma}

\begin{proof}
	Let $I_{j-1}=(\alpha_1,\ldots,\alpha_{j-1})$. Recall that for a monomial ideal
	$J=(u_1,\ldots,u_t)$ and a monomial $v$,
	\[
	J:v=\left(\frac{u_1}{\gcd(u_1,v)},\ldots,
	\frac{u_t}{\gcd(u_t,v)}\right);
	\]
	see also \cite[Proposition 1.2.2]{HH}. We apply this formula with $J=I_{j-1}$ and $v=\alpha_j$, and remove
	redundant monomials from the resulting generating set.
	
	1. If $\alpha_j=e_{a,1}e_{b,2}$ for some $1\leq a<b\leq m-1$, then the
	preceding generators $e_{i,1}e_{b,2}$, $1\leq i<a$, give the quotients
	$e_{i,1}$, while $e_{a,1}e_{k,2}$, $a<k<b$, give the quotients $e_{k,2}$.
	Every other quotient is divisible by one of these monomials. Hence
	\[
	I_{j-1}:\alpha_j
	=(e_{1,1},\ldots,e_{a-1,1},e_{a+1,2},\ldots,e_{b-1,2}).
	\]
	
	2. If $\alpha_j=g_1f_2e_{a,1}$ for some $1\leq a\leq m-1$, then
	\[
	I_{j-1}:\alpha_j
	=(e_{1,1},\ldots,e_{a-1,1},e_{a+1,2},\ldots,e_{m-1,2}).
	\]
	The displayed generators arise from $g_1f_2e_{i,1}$ with $i<a$ and
	$e_{a,1}e_{k,2}$ with $k>a$, respectively. For every other generator
	$e_{i,1}e_{k,2}\in\mathcal M_1$, either $i<a$ or $k>a$, so the
	corresponding quotient is divisible by one of the displayed generators.
	
	3. If $\alpha_j=f_1f_2g_1g_2$, then
	\[
	I_{j-1}:\alpha_j=(e_{1,1},\ldots,e_{m-1,1}).
	\]
	Indeed, the generators $g_1f_2e_{a,1}$ give the quotients $e_{a,1}$, and every other quotient is divisible by one of them.
	
	4. If $\alpha_j=g_2f_1e_{a,2}$ for some $1\leq a\leq m-1$, then
	\[
	I_{j-1}:\alpha_j
	=(g_1f_2,e_{1,1},\ldots,e_{a-1,1},e_{a+1,2},\ldots,e_{m-1,2}).
	\]
	The displayed generators arise from $f_1f_2g_1g_2$,
	$e_{i,1}e_{a,2}$ with $i<a$, and $g_2f_1e_{t,2}$ with $t>a$,
	respectively. The quotients arising from the generators in
	$\mathcal M_3$ are divisible by $g_1f_2$. For any other generator
	$e_{i,1}e_{k,2}\in\mathcal M_1$, either $i<a$ or $k>a$, so the
	corresponding quotient is divisible by one of the displayed generators.
	
	In each case, the listed generators have pairwise disjoint supports, and hence
	form a regular sequence. Counting the displayed generators in the four cases gives
	$k_j=b-2$, $m-2$, $m-1$, and $m-1$, respectively.
\end{proof}

\begin{Lemma} \label{Gmupper}
Let $G$ be the graph as described above. Then there exists a monomial order $<$ on
$\mathbb{K}[E(G)]$ such that, letting $I:=\mathrm{in}_{<}(I_{G})$, we have
\[
  \beta_i(\mathbb{K}[G])
  \le \beta_i\bigl(\mathbb{K}[E(G)]/I\bigr)
  \le i\binom{m+1}{i+1}
  \quad \text{for all } i\ge 1.
\]
\end{Lemma}

\begin{proof}
Let $R=\mathbb{K}[E(G)]$, and let
$\alpha_1 \prec \cdots \prec \alpha_r$ denote  the ordering of the minimal generators of $I$ given prior to Lemma~\ref{Gmin}.
For $s=1,\ldots,r$, set
\[
  I_s := (\alpha_1,\ldots,\alpha_s)
  \quad\text{and}\quad
  d_s := \deg(\alpha_s),
\]
and for $s\ge 2$ define
\[
  L_s := (\alpha_1,\ldots,\alpha_{s-1}):\alpha_s.
\]
By Lemma~\ref{Gmin}, each $L_s$ is generated by a regular sequence of length $k_s$.

Since $L_s$ is generated by a regular sequence of length $k_s$, the minimal free
resolution of $R/L_s$ is the Koszul complex on this regular sequence. Hence, for all $i\ge 0$,
\[
  \beta_i(R/L_s) = \binom{k_s}{i}.
\]

For each $s=2,\ldots,r$, consider the short exact sequence
\[
0 \longrightarrow R/L_s(-d_s)
   \xrightarrow{\ \times \alpha_s\ }
   R/I_{s-1}
   \longrightarrow R/I_s
   \longrightarrow 0,
\]
where the first map is multiplication by $\alpha_s$.
Applying the mapping cone construction (see, for example, \cite[Section~27]{P}),
we obtain, for all $i\ge 1$,
\[
\beta_i(R/I_s)
\;\leq\; \beta_{i-1}(R/L_s) + \beta_i(R/I_{s-1}).
\]
Here and in what follows we work with total Betti numbers, that is, we ignore the internal grading.

Iterating this inequality on $s$ and using that $I_1=(\alpha_1)$ is principal,
one obtains
\[
  \beta_i(R/I)
  \;\leq\; \sum_{s=2}^r \beta_{i-1}(R/L_s) + \beta_i(R/I_1)
  \quad\text{for all } i\ge 1.
\]
Since $I_1$ is principal, the one-element Koszul resolution gives $\beta_i(R/I_1)=\binom{1}{i}$.
Substituting $\beta_{i-1}(R/L_s)=\binom{k_s}{i-1}$, we get
\[
  \beta_i(R/I)
  \;\leq\; \sum_{s=2}^r \binom{k_s}{i-1} + \binom{1}{i}.
\]

From Lemma~\ref{Gmin}, we obtain the following combinatorial fact:
for every integer $k$ with $1\leq k\leq m-1$, the value $k$ occurs exactly
$k+1$ times among
\[
k_2,k_3,\ldots,k_r.
\]
Indeed, when $m=1$ there are no such $k$, and when $m=2$ we have
$(k_2,k_3)=(1,1)$. Assume that $m\geq 3$. For $1\leq k\leq m-3$,
the value $k$ corresponds precisely to the generators
$e_{a,1}e_{k+2,2}$ with $1\leq a\leq k+1$, and hence occurs $k+1$ times.
The value $m-2$ occurs $m-1=(m-2)+1$ times, corresponding to the generators
$g_1f_2e_{a,1}$ with $1\leq a\leq m-1$. Finally, the value $m-1$ occurs
$m=(m-1)+1$ times, corresponding to $f_1f_2g_1g_2$ and the $m-1$
generators $g_2f_1e_{a,2}$.

Therefore,
\[
  \sum_{s=2}^r \binom{k_s}{i-1}
    = \sum_{k=1}^{m-1} (k+1)\binom{k}{\,i-1\,}.
\]
Combining this with the above inequality for $\beta_i(R/I)$, we obtain
\begin{align*}
  \beta_{i}(R/I)
  &\leq \sum_{k=1}^{m-1} (k+1)\binom{k}{\,i-1\,} + \binom{1}{i} \\[2pt]
  &= i\sum_{k=1}^{m-1} \binom{k+1}{i} + \binom{1}{i} \\[2pt]
  &= i\sum_{k=0}^{m-1} \binom{k+1}{i}
   \;=\; i\binom{m+1}{i+1}.
\end{align*}
Here, the first equality uses the binomial identity
$(k+1)\binom{k}{i-1} = i\binom{k+1}{i}.$

Finally, the inequalities for $\beta_i\bigl(\mathbb{K}[E(G)]/I\bigr)$ and $\beta_i(\mathbb{K}[G])$
follow from Lemma~\ref{total}, applied to $I_G$ and its initial ideal
$I=\mathrm{in}_{<}(I_G)$.
\end{proof}

\begin{Remark}\em
Lemma~\ref{Gmin} shows that the above ordering of the minimal generators of
$I=\mathrm{in}_{<}(I_G)$ gives $I$ the property called having \emph{regular quotients}
in \cite{LZ, M20}.
\end{Remark}

\subsection{Conventions and notation}
We first introduce the necessary conventions and notation, and clarify
key notions used in the induced-subgraph approach. These preliminaries
will also be utilized in subsequent sections.

\begin{Conventions}\label{conventions}\em
Let $G$ be a simple graph. We write $\RR^{V(G)}$ for the real vector space of all functions from $V(G)$ to $\RR$, with its standard basis indexed by the vertices of $G$. If the vertices are labeled by letters or subscripted letters, such as $u, u_{i,j}, v,  v_{i,j}, w, w_{i,j} $, and so forth, we establish a correspondence between these labels and the corresponding unit vectors in $\RR^{V(G)}$ as follows: $$u, u_{i,j}, v,  v_{i,j}, w, w_{i,j}$$ correspond to $$\ub, \ub_{i,j}, \vb, \vb_{i,j}, \wb, \wb_{i,j}$$ respectively. With this notation, we can express $\RR^{V(G)}$ as
the set\begin{align*}
 \left\{ \sum_{i,j} a_{i,j}\ub_{i,j} + a\ub + \sum_{i,j} b_{i,j}\vb_{i,j} + b\vb + \sum_{i,j} c_{i,j}\wb_{i,j} + c\wb \mid a_{i,j}, a, b_{i,j}, b, c_{i,j}, c \in \RR\right\}.
\end{align*}
By way of example,  $\deg(uv)=\ub + \vb$, and $\deg(u^2wu_{1,2})=2\ub + \wb + \ub_{1,2}$, and so on.

We identify a monomial with its multidegree. For example, we will  write $\beta_{i,uv}(\mathbb{K}[G])$ for $\beta_{i, \ub+\vb}(\mathbb{K}[G])$.
 Under this identification, we regard $\mathbb{Z}^{V(H)}$ as a subset of $\mathbb{Z}^{V(G)}$
whenever $H$ is an induced subgraph of $G$.
\end{Conventions}

\begin{Notation}\em
	Let $H$ be a simple graph satisfying the odd cycle condition. For every edge
	$e=\{u,v\}$ of $H$, let $\mathbf{v}_e=\mathbf{u}+\mathbf{v}$. We define
	$D_H$ to be the multidegree
	$\sum_{e\in E(H)}\mathbf{v}_e\in \mathbb{Z}_{\geq 0}^{V(H)}$.
	Since we identify a monomial with its multidegree, $D_H$ can alternatively
	be interpreted as the product of all edges of $H$. Here, an edge is considered
	as the product of its two vertices. Furthermore, we employ $\Theta_H$ to
	represent the product of all the vertices of $H$.
\end{Notation}
If $H$ is a simple graph in which every vertex has degree at least $2$, then it is straightforward to see that $\Theta_H^{2}$ divides $D_H$.
 Moreover, a simple graph $H$ is a top-Betti graph if and only if $\Theta_H$ divides every monomial in $\mathcal{N}_H$.

Let $G$ be a simple graph, and denote the set of all monomials in the edge ring $\mathbb{K}[G]$ of $G$ by $\mathrm{Mon}(\mathbb{K}[G])$.

\begin{Definition}\em
Let $G$ be a simple graph with vertex set $V(G)=\{u_1,\ldots,u_n\}$ and edge set $E(G)$.
For each edge $e=\{u_i,u_j\}\in E(G)$, define
\[
\mathbf{v}_{e} := \mathbf{u}_{i}+\mathbf{u}_{j},
\]
where $\mathbf{u}_{i}$ denotes the $i$th unit vector in $\mathbb{R}^{V(G)}$.
The \emph{edge cone} of $G$, denoted by $\mathbb{R}_{+}(G)$, is the cone in $\mathbb{R}^{V(G)}$
spanned by the vectors $\{\mathbf{v}_{e}\mid e\in E(G)\}$; explicitly,
\[
\mathbb{R}_{+}(G)
:= \left\{ \sum_{e\in E(G)} a_e \mathbf{v}_e \;\bigg|\;
a_e \in \mathbb{R}_{\ge 0} \text{ for all } e\in E(G) \right\}.
\]
We denote by $\mathrm{relint}\bigl(\mathbb{R}_{+}(G)\bigr)$ the relative interior of the cone $\mathbb{R}_{+}(G)$.
\end{Definition}

  It was established by \cite[Corollary 2.3]{OH} that if $G$ satisfies the odd cycle condition, then the edge ring $\mathbb{K}[G]$ is normal. Let  $\mathcal{S}(G)$ denote the semigroup composed of all $h\in \mathbb{Z}^{V(G)}$ with $X^{h}\in \mathbb{K}[G].$ According to \cite[Theorem~6.31]{BG}, the ideal of $\mathbb{K}[G]$
 generated by all monomials $X^{h}$ with
 $h \in \mathrm{relint}(\mathbb{R}_{+}(G)) \cap \mathcal{S}(G)$
 constitutes the $\mathbb{Z}^{V(G)}$-graded canonical module of the
 $\mathbb{Z}^{V(G)}$-graded $^*$local algebra $\mathbb{K}[G]$, which we
 denote by $\omega_{\mathbb{K}[G]}$.
 Set $S=\mathbb{K}[E(G)]$ and let
 $p=\operatorname{pdim}(\mathbb{K}[G])$. Since the sum of the
 multidegrees of the variables of $S$ is $D_G$, its multigraded
 canonical module is $\omega_S\cong S[-D_G]$. Since
 $\mathbb{K}[G]$ is Cohen--Macaulay, its codimension in $S$ is $p$.
 Applying \cite[Equation~6.6]{BG}, together with
 \[
 \operatorname{Hom}_S(S[-h],\omega_S)
 \cong S[-(D_G-h)],
 \]
 yields the following.

\begin{Lemma}\label{canonical} Let $G$ be a  simple graph  satisfying the odd cycle condition and denote by $p$ the projective dimension of $\mathbb{K}[G]$.  Then $$\beta_{p-i,h}(\mathbb{K}[G])=\beta_{i,D_G-h}(\omega_{\mathbb{K}[G]})\qquad  \mbox{ for all }\quad 0\leq i\leq p\mbox{ and for all }h\in \mathbb{Z}_{\geq 0}^{V(G)}.$$
\end{Lemma}

\subsection{Returning to $\mathbf{G}_m$}
\begin{Proposition}\label{prop:two-ear-top-betti-types}
Let $G$ be a two-ear graph $\mathbf{G}_m$. Then every top-Betti induced subgraph $H$ of $G$ is one of the following three types:
\begin{itemize}
  \item a two-ear graph;
  \item a one-ear graph;
  \item a complete bipartite graph of the form $\mathbf{K}_{2,d}$ for some $d\geq 2$.
\end{itemize}
\end{Proposition}

\begin{proof}
By Proposition~\ref{Betti}, every two-ear graph, every one-ear graph
of type $r$ with $r\geq 2$, and every complete bipartite graph
$\mathbf{K}_{2,d}$ with $d\geq 2$ is a top-Betti graph. We next prove
that a two-ear graph has no other top-Betti induced subgraphs.

Label the vertices and edges of $G$ as in Figure~1, and let $H$ be
a top-Betti induced subgraph of $G$. By Corollary~\ref{top-Betti at least two},
every vertex of $H$ has degree at least $2$. In particular,
if $x\in V(H)$ and $\deg_G(x)=2$, then all neighbors of $x$ in $G$ must also
belong to $V(H)$, that is,
\[
N_G(\{x\}) \subseteq V(H).
\]

Set $F:=\{u_1,\ldots,u_{m-1}\}$. We distinguish two cases depending on whether
$V(H)\cap F$ is empty or not.

\medskip
\noindent\textbf{Case 1.}
$V(H)\cap F\neq\emptyset$.

Since each vertex $u_i\in F$ satisfies $\deg_G(u_i)=2$ and $N_G(\{u_i\})=\{v_1,v_2\}$,
the above observation implies $v_1,v_2\in V(H)$.
If further
\[
(V(H)\cap F)\sqcup\{v_1,v_2\}=V(H),
\]
then $|V(H)\cap F|\ge 2$ holds; for if this is not the case, then $\deg_H(v_1)=1$, which contradicts the fact that $\deg_H(v_1)\ge 2$.
It follows that $H$ is the complete bipartite graph with bipartition $\{v_1,v_2\}$ and $V(H)\cap F$, i.e., a graph of the form $\mathbf{K}_{2,d}$ where $d = |V(H)\cap F|\geq 2$.

Now suppose $(V(H)\cap F)\sqcup\{v_1,v_2\}$ is a proper subset of $V(H)$.
This is equivalent to the condition
\[
\bigl|V(H)\cap\{x_1,x_2,y_1,y_2\}\bigr|\ge 1.
\]
The case
\[
\bigl|V(H)\cap\{x_1,x_2,y_1,y_2\}\bigr|=1
\]
is impossible, since it would force $H$ to contain a vertex of degree $1$.
Note that $\deg_G(x_2)=\deg_G(y_2)=2$, and that $\{x_1,x_2\}$ and $\{y_1,y_2\}$ are edges of $G$.
Thus, if $x_2$ (resp.\ $y_2$) lies in $V(H)$, then $x_1$ (resp.\ $y_1$) must also lie in $V(H)$.
Consequently, the intersection
\[
V(H)\cap\{x_1,x_2,y_1,y_2\}
\]
can only take one of the following forms:
\[
\begin{aligned}
&\{x_1,x_2\},\quad \{y_1,y_2\},\quad \{x_1,y_1\},\\
&\{x_1,y_1,y_2\},\quad \{x_1,x_2,y_1\},\quad \{x_1,x_2,y_1,y_2\}.
\end{aligned}
\]

If $V(H)\cap\{x_1,x_2,y_1,y_2\}$ is one of
\[
\{x_1,y_1,y_2\},\quad \{x_1,x_2,y_1\},\quad \{x_1,x_2,y_1,y_2\},
\]
then $H$ is either a one-ear graph or a two-ear graph.

Next, we show that the remaining three possibilities
\[
\{x_1,x_2\},\ \{y_1,y_2\},\ \{x_1,y_1\}
\]
cannot occur.
Indeed, if
\[
V(H)\cap\{x_1,x_2,y_1,y_2\}\in\bigl\{\{x_1,x_2\},\ \{y_1,y_2\},\ \{x_1,y_1\}\bigr\},
\]
then the induced subgraph $H$ is non-bipartite.
Moreover, for every $x\in V(H)\cap\{x_1,x_2,y_1,y_2\}$, we have $\deg_H(x)=2$ and $H\setminus x$ is bipartite.
Hence
\[
c_0(H\setminus x)=1=\deg_H(x)-1,
\]
which contradicts Proposition~\ref{Betti} (applied to the non-bipartite graph $H$).
Therefore, in Case~1, $H$ must be either a complete bipartite graph of the form
$\mathbf{K}_{2,d}$, a one-ear graph, or a two-ear graph.

\medskip
\noindent\textbf{Case 2.} $V(H)\cap F=\emptyset$.
In this case, $H$ is an induced subgraph of the graph obtained from $G$
by deleting all vertices of $F$, which is $\mathbf{G}_1$. The graph $\mathbf{G}_1$ consists of the two triangles on $\{v_1,x_1,x_2\}$ and $\{v_2,y_1,y_2\}$, joined by the edge $\{x_1,y_1\}$. We claim that a proper induced subgraph of $\mathbf{G}_1$ whose vertices all have degree at least $2$ must be one of these two triangles.

Indeed, let $L$ be such a proper induced subgraph. If $L$ contains $v_1$ or
$x_2$, then the degree condition forces it to contain the whole triangle
$\{v_1,x_1,x_2\}$. If $x_1\in V(L)$ but $v_1,x_2\notin V(L)$, then
$\deg_L(x_1)\leq1$. Hence, if $L$ meets the triangle
$\{v_1,x_1,x_2\}$, it contains the whole triangle. The same argument applies
to $\{v_2,y_1,y_2\}$. Therefore, $L$ contains at least one of the two
triangles, while it cannot contain both, since otherwise inducedness would
give $L=\mathbf{G}_1$. Thus $L$ is one of the two triangles.

Since a triangle is not a top-Betti graph by Proposition~\ref{Betti}, $H$
cannot be a proper induced subgraph of $\mathbf{G}_1$. Hence
$H=\mathbf{G}_1$, which is a two-ear graph.

Combining the two cases, the proposition follows.
\end{proof}

Next, we determine the elements of the top support for the three types of graphs by considering the canonical module of their respective edge rings. For the definition of top support, we refer to Definition~\ref{topsupport}. By Lemma~\ref{canonical}, every nonzero top multigraded Betti number occurring here is equal to 1. Hence the top support completely determines the top multigraded
Betti numbers. The most complex scenario is when $H$ is a two-ear graph, for which we will provide a detailed proof. When $H$ is a one-ear graph or a complete bipartite graph $\mathbf{K}_{2,d}$, the proofs are similar with only minor differences, so we only outline the proofs for those cases.

\subsection{The top support of $\mathbf{G}_m$}
Notice that when $m=1$, the two-ear graph $\mathbf{G}_1$ is isomorphic to the compact graph $\mathbf{B}^{0}_{(1):(1)}$, which we discuss in the next section. Its top support is as shown in Lemma~\ref{5.3}.
Therefore, in this subsection, we assume $m \geq 2$ and, for brevity, use $G$ to represent $\mathbf{G}_m$.

It is easy to see that  $|V(G)|=m+5$. It follows from Conventions~\ref{conventions} that
\begin{align*}
	\RR^{V(G)}=\{\sum\limits_{i=1}^{m-1}{a_{i}\ub_{i}}+
	\sum\limits_{i=1}^{2}({b_{i}\vb_{i}+c_{i}\xb_{i}+d_{i}\yb_{i}})
	\mid a_{i}, b_{i},c_{i},d_{i}\in \RR \mbox{ for all } i \},
\end{align*}
where $\ub_{i},\vb_{i},\xb_{i},\yb_{i}$ are the unit vectors of $\RR^{V(G)}$,
each vector $\ub_{i}$ (resp. $\vb_{i},\xb_{i},\yb_{i}$) corresponds to vertex $u_{i}$ (resp. $v_{i},x_{i},y_{i}$) for all $i$.

Recall that an integral vector in $\mathrm{relint}(\RR_+(G))\cap \mathcal{S}(G)$ is  {\it minimal} if it cannot be written as the sum of a vector in $\mathrm{relint}(\RR_+(G))\cap \mathcal{S}(G)$ and a nonzero vector of  $\RR_+(G)\cap \mathcal{S}(G)$. It follows that a minimal vector of $\mathrm{relint}(\RR_+(G))\cap \mathcal{S}(G)$ corresponds to a minimal generator of the canonical module of $\mathbb{K}[G]$.

To find the minimal vectors in $\mathrm{relint}(\RR_+(G))\cap \mathcal{S}(G)$, we need to express the edge cone $\RR_{+}(G)$ in terms of linear inequalities,
or equivalently, through its facets.

\begin{Definition}\em\label{def-regular-fundamental}
Let $G$ be a connected non-bipartite simple graph.
\begin{enumerate}
  \item A vertex $j$ of $G$ is called \emph{regular} if each connected component of
  $G \setminus j$ contains an odd cycle.
  \item An independent set $T\subseteq V(G)$ is called \emph{fundamental} if it satisfies:
  \begin{itemize}
    \item[(a)] the bipartite graph on the vertex set $T \sqcup N(T)$ with edge set
    \[
     \bigl\{\{j,k\}\in E(G)\mid j\in T\bigr\}
     \]
    is connected; and
    \item[(b)] either $T \sqcup N(T)=V(G)$, or each connected component of
    $G\bigl[V(G) \setminus (T \sqcup N(T))\bigr]$ contains an odd cycle.
  \end{itemize}
\end{enumerate}
\end{Definition}

We hereby record  a result from  \cite[Theorem 10.7.6]{V} as follows:

 \begin{Lemma} \label{4.12} Let $G$ be a non-bipartite connected simple graph on $V(G)=\{1,\ldots,n\}$. Then
 $\RR_{+}(G)$ consists of all elements $(x_1, \ldots, x_n) \in \RR^n$ satisfying the following inequalities:
\begin{equation}\label{eq:ineq2}\tag{$\Delta$}
\begin{split}
x_i &\geq 0 \;\;\text{ for any regular vertex }i; \\
\sum\limits_{i \in N(T)}x_i &\geq \sum\limits_{i \in T}x_i \;\;\text{ for any fundamental set }T.
\end{split}
\end{equation}
\end{Lemma}

\begin{Lemma}\label{fundamental-T}
Let $G$ be the graph $\mathbf{G}_m$ described above. Then all vertices of $G$ are regular. Moreover,
an independent subset $T$ of $V(G)$ is fundamental if and only if $T$ is one of the following:
\begin{enumerate}
  \begin{multicols}{2}
    \item[$(i)$] $\{v_1\};$
    \item[$(ii)$] $\{v_2\};$
    \item[$(iii)$] $\{x_2\};$
    \item[$(iv)$] $\{y_2\};$
    \item[$(v)$] $\{v_1,v_2\};$
    \item[$(vi)$] $\{v_1,y_1\};$
    \item[$(vii)$] $\{v_2,x_1\};$
    \item[$(viii)$] $\{u_1,u_2,\ldots,u_{m-1},x_1,y_2\};$
    \item[$(ix)$] $\{u_1,u_2,\ldots,u_{m-1},x_2,y_1\};$
    \item[$(x)$] $\{u_1,u_2,\ldots,u_{m-1},x_2,y_2\}.$
  \end{multicols}
  \end{enumerate}

\begin{proof}
By Definition~\ref{def-regular-fundamental}(1), a vertex $j$ of $G$ is regular if every connected component of $G\setminus j$ contains an odd cycle. A straightforward verification shows that this criterion holds for every vertex of the two-ear graph $G$, whence all vertices of $G$ are regular.

To establish the second assertion, we classify the fundamental independent sets of $G$. By Definition~\ref{def-regular-fundamental}(2), an independent set $T\subseteq V(G)$ is fundamental if and only if it satisfies conditions~(a) and~(b) therein. We proceed to characterize all such sets explicitly.

For any $z\in V(G)$, the singleton $\{z\}$ is independent and satisfies (a), since $G$ is connected.
Moreover, by the definition of $\mathbf{G}_m$, $z$ is not adjacent to all the other vertices of $G$, and hence
\[
\{z\}\sqcup N(\{z\}) \neq V(G).
\]
Thus $\{z\}$ is fundamental if and only if
\[
  G\bigl[V(G) \setminus (\{z\} \sqcup N(\{z\}))\bigr]
\]
contains an odd cycle. A direct verification shows that this happens precisely
for the sets in \emph{(i)}-\emph{(iv)}.

Now let $T$ be a fundamental independent set with $|T|>1$. If $v_1\in T$, then the only
non-neighbors of $v_1$ are $y_1,y_2,v_2$, which form a triangle in $G$. Hence
\[
T\in\bigl\{\{v_1,y_1\},\ \{v_1,y_2\},\ \{v_1,v_2\}\bigr\}.
\]
One checks that $\{v_1,y_1\}$ and $\{v_1,v_2\}$ satisfy (a) and (b). On the other hand,
\[
N(\{y_2\})=\{y_1,v_2\},\qquad
N(\{v_1\})=\{u_1,\ldots,u_{m-1},x_1,x_2\}
\]
are disjoint. Thus, for $T=\{v_1,y_2\}$, the bipartite graph on $T\sqcup N(T)$ is disconnected,
so (a) fails. Similarly, if $v_2\in T$, then $T$ must be either $\{v_2,x_1\}$ or $\{v_1,v_2\}$.
Hence, when $|T|>1$ and $v_1$ or $v_2$ belongs to $T$, the only possibilities are
\emph{(v)}-\emph{(vii)}.

We now assume that $v_1,v_2\notin T$. If
\[
T\cap\{u_1,\ldots,u_{m-1}\}=\emptyset,
\]
then $T$ can only be
\[
\{x_2,y_1\},\quad \{x_1,y_2\},\quad \{x_2,y_2\}.
\]
If $T=\{x_2,y_1\}$ or $T=\{x_1,y_2\}$, then $T\sqcup N(T)\neq V(G)$ and the induced subgraph
\[
G\bigl[V(G) \setminus (T \sqcup N(T))\bigr]
\]
has vertex set $\{u_1,\ldots,u_{m-1}\}$ and no edges, so it contains no odd cycle, contradicting (b).
If $T=\{x_2,y_2\}$, the neighborhoods of $x_2$ and $y_2$ are disjoint, so the bipartite graph on
$T\sqcup N(T)$ is disconnected, contradicting (a). Thus
\[
T\cap\{u_1,\ldots,u_{m-1}\}\neq\emptyset.
\]

In particular, $v_1,v_2\notin V(G)\setminus(T\sqcup N(T))$. Every odd cycle of $G$ passes through
$v_1$ or $v_2$, hence the induced subgraph
\[
G\bigl[V(G) \setminus (T \sqcup N(T))\bigr]
\]
is bipartite. By condition~(b), this forces $T\sqcup N(T)=V(G)$. Since each $u_i$ is adjacent only to $v_1$
and $v_2$, we must have
\[
\{u_1,\ldots,u_{m-1}\}\subseteq T,
\]
for otherwise some $u_i\notin T$ would also satisfy $u_i\notin N(T)$ (because $v_1,v_2\notin T$),
contradicting $T\sqcup N(T)=V(G)$. Consequently, $T$ can only be
\[
\{u_1,\ldots,u_{m-1},x_1,y_2\},\quad
\{u_1,\ldots,u_{m-1},x_2,y_1\},\quad
\{u_1,\ldots,u_{m-1},x_2,y_2\},
\]
i.e., the cases \emph{(viii)}-\emph{(x)}. A direct verification shows that they all satisfy (a),
so the proof is complete.
\end{proof}

\end{Lemma}

By Lemma~\ref{fundamental-T} and Lemma~\ref{4.12}, it follows  that a vector in $\RR^{V(G)}$ of the form:
$$\sum\limits_{i=1}^{m-1}{a_{i}\ub_{i}}+ \sum\limits_{i=1}^{2}({b_{i}\vb_{i}+c_{i}\xb_{i}+d_{i}\yb_{i}})$$ belongs to $ \RR_+(G)$ if and only if it satisfies the following inequalities:

\begin{center}
	\begin{tabular}{@{}p{0.48\textwidth}p{0.48\textwidth}@{}}
		(1) $a_i\geq 0$ for any $1\leq i\leq m-1$;
		&
		(2) $b_i,c_i,d_i\geq 0$ for any $1\leq i\leq 2$;
		\\[1.5mm]
		
		(3) $\displaystyle\sum_{i=1}^{m-1}a_i+c_1+c_2\geq b_1$;
		&
		(4) $\displaystyle\sum_{i=1}^{m-1}a_i+d_1+d_2\geq b_2$;
		\\[1.5mm]
		
		(5) $c_1+b_1\geq c_2$;
		&
		(6) $d_1+b_2\geq d_2$;
		\\[1.5mm]
		
		\multicolumn{2}{@{}l@{}}{
			(7) $\displaystyle\sum_{i=1}^{m-1}a_i+c_1+c_2+d_1+d_2\geq b_1+b_2$;
		}
		\\[1.5mm]
		
		\multicolumn{2}{@{}l@{}}{
			(8) $\displaystyle\sum_{i=1}^{m-1}a_i+b_2+c_1+c_2+d_2\geq b_1+d_1$;
		}
		\\[1.5mm]
		
		\multicolumn{2}{@{}l@{}}{
			(9) $\displaystyle\sum_{i=1}^{m-1}a_i+b_1+c_2+d_1+d_2\geq b_2+c_1$;
		}
		\\[1.5mm]
		
		\multicolumn{2}{@{}l@{}}{
			(10) $b_1+b_2+c_2+d_1\geq
			\displaystyle\sum_{i=1}^{m-1}a_i+c_1+d_2$;
		}
		\\[1.5mm]
		
		\multicolumn{2}{@{}l@{}}{
			(11) $b_1+b_2+c_1+d_2\geq
			\displaystyle\sum_{i=1}^{m-1}a_i+c_2+d_1$;
		}
		\\[1.5mm]
		
		\multicolumn{2}{@{}l@{}}{
			(12) $b_1+b_2+c_1+d_1\geq
			\displaystyle\sum_{i=1}^{m-1}a_i+c_2+d_2$.
		}
	\end{tabular}
\end{center}
We now construct $m$ integral vectors in $\RR^{V(G)}$ and then show that they are minimal vectors of $\mathrm{relint}(\RR_+(G))\cap \mathcal{S}(G)$.
The construction is as follows:

For $\ell=1,\ldots,m$, let \begin{equation*}
\alpha_{\ell}:=\sum_{i=1}^{m-1}\ub_{i}+\sum_{i=1}^{2}\bigl(\xb_{i}+\yb_{i}\bigr)+\ell\vb_{1}+(m+1-\ell)\vb_{2}.
\end{equation*}

It is straightforward to check that each $\alpha_{\ell}$ satisfies all twelve of the required inequalities strictly. While a full verification of all inequalities would be redundant, we present a detailed check for inequality (7) as a representative example. For any $\alpha_{\ell}$, we have
\[
a_i=1\ (1\le i\le m-1),\quad c_1=c_2=d_1=d_2=1,\quad b_1=\ell,\quad b_2=m+1-\ell.
\]
Substituting these values into inequality (7) yields
\[
\sum_{i=1}^{m-1}a_i+c_1+c_2+d_1+d_2
=(m-1)+4=m+3
> m+1=\ell+(m+1-\ell)=b_1+b_2.
\]
Thus the inequality (7) holds strictly.

It is also obvious that the sum of all components of $\alpha_{\ell}$ is even, and according to \cite[Proposition 3.4]{V}, we have $\alpha_{\ell} \in \mathrm{relint}(\RR_+(G))\cap \mathcal{S}(G)$.

Next, we show that each $\alpha_{\ell}$ is a minimal vector in
$\mathrm{relint}(\RR_+(G))\cap \mathcal{S}(G)$. Suppose, to the contrary, that
\[
\alpha_{\ell}=\alpha' + \alpha''
\]
for some $\alpha' \in \mathrm{relint}(\RR_+(G))\cap \mathcal{S}(G)$ and
$\alpha'' \in \RR_+(G)\cap \mathcal{S}(G) \setminus \{\mathbf{0}\}$.
We write
\[
\alpha'=\sum_{i=1}^{m-1} a'_{i}\ub_{i}
        + \sum_{i=1}^{2}\bigl(b'_{i}\vb_{i}+c'_{i}\xb_{i}+d'_{i}\yb_{i}\bigr),
\qquad
\alpha''=\sum_{i=1}^{m-1} a''_{i}\ub_{i}
         + \sum_{i=1}^{2}\bigl(b''_{i}\vb_{i}+c''_{i}\xb_{i}+d''_{i}\yb_{i}\bigr).
\]

Note that
\[
\mathrm{relint}(\RR_+(G))\cap \mathcal{S}(G)
 \subseteq \RR_+(G)\cap \mathcal{S}(G) \subseteq \mathbb{Z}^{V(G)}.
\]
In view of inequalities (1) and (2), we see that
$a'_{i}, c'_{i}, d'_{i}\geq 1$ and $a''_{i}, c''_{i}, d''_{i}\geq 0$
for all $i$. Hence $a'_{i}=c'_{i}=d'_{i}=1$ and
$a''_{i}=c''_{i}=d''_{i}=0$ for all $i$.
Because of (7), we have $b''_{1} + b''_{2} \leq 0$.
But from (2) we know $b''_1, b''_2\geq 0$, so we obtain $b''_1=b''_2=0$.

Hence $\alpha''=\mathbf{0}$, a contradiction. Therefore
$\alpha_{\ell}$ is a minimal vector in
$\mathrm{relint}(\RR_+(G))\cap \mathcal{S}(G)$ for $\ell=1,\ldots,m$.

Actually, $\alpha_{\ell}$ for $\ell=1,\ldots,m$ are all the minimal vectors in
$\mathrm{relint}(\RR_+(G))\cap \mathcal{S}(G)$, as shown below.

\begin{Lemma}\label{GmNH}
Let $m \geq 1$. Let $G$ denote the graph $\mathbf{G}_m$ with the labeling as shown in Figure~\ref{fig3}. Then we have
\[
\mathrm{pdim}(\mathbb{K}[G]) \;=\;
\mathrm{pdim}\bigl(\mathbb{K}[E(G)] / \mathrm{in}_{<}(I_G)\bigr)
\;=\; m.
\]
Moreover,
\[
\mathcal{N}_G
= \bigl\{
\Theta_G \, x_1 y_1 \, v_1^{\,m-\ell} v_2^{\,\ell-1}
\;\bigm|\;
\ell = 1,\ldots,m
\bigr\}.
\]
\end{Lemma}

\begin{proof}
First suppose that $m=1$. In this case $G$ is the graph obtained by joining two vertex-disjoint triangles by a single edge. By Lemma~\ref{pd} we have
\[
\mathrm{pdim}(\mathbb{K}[G]) = 1.
\]
The toric ideal $I_G$ is generated by the single binomial
$f_1 f_2 g_1 g_2 - f_3 g_3 w^2$, and the multidegree of this binomial is
\[
v_1 x_2 x_1 y_1 y_2 v_2 x_1 y_1
  \;=\; \Theta_G x_1 y_1.
\]
Hence
\[
\beta_{1,h}(\mathbb{K}[G]) = \beta_{0,h}(I_G) \neq 0
\quad\Longleftrightarrow\quad
X^h = \Theta_G x_1 y_1,
\]
which yields the desired description of $\mathcal{N}_G$ when $m=1$.

Now assume that $m \geq 2$. The first assertion
\[
\mathrm{pdim}(\mathbb{K}[G])
= \mathrm{pdim}\bigl(\mathbb{K}[E(G)] / \mathrm{in}_{<}(I_G)\bigr)
= m
\]
follows from Lemma~\ref{pd} together with \cite[Corollary~2.7]{CV}. For the second assertion, note that since
\(\mathrm{pdim}(\mathbb{K}[G]) = m\), the top total Betti number of $\mathbb{K}[G]$ is at most $m$ by Lemma~\ref{Gmupper}. Consequently, the cardinality of the set \(\mathcal{N}_G\) satisfies
$|\mathcal{N}_G| \leq m.$
On the other hand, from the previous discussion we already know that
$|\mathcal{N}_G| \geq m.$
Moreover, we have
\[
D_G = \Theta_G^2 x_1 y_1 (v_1 v_2)^{m-1}.
\]
Applying Lemma~\ref{canonical}, we conclude that
\[
\mathcal{N}_G
= \bigl\{
\Theta_G x_1 y_1 v_1^{\,m-\ell} v_2^{\,\ell-1}
\;\bigm|\;
\ell = 1,\ldots,m
\bigr\},
\]
as desired.
\end{proof}

A direct computation shows that every monomial in $\mathcal{N}_{\mathbf{G}_m}$ has degree $2m+6$.
For example, by \cite[Corollary~2.17]{HHO} the regularity of a Cohen--Macaulay module is determined by its top Betti numbers. Combining this result with the relation between the multigrading and the standard grading described in the Introduction, we obtain the following.
\begin{Corollary}
The regularity of the edge ring $\mathbb{K}[\mathbf{G}_m]$ is equal to $3$.
\end{Corollary}

\subsection{The one-ear graph $\mathbf{P}_m$}

 Let $m \geq 2$ be an integer. Let $G'$ (resp. $G''$) denote the induced subgraph $\mathbf{G}_m\setminus x_2$ (resp. $\mathbf{G}_m\setminus y_2 $) of $\mathbf{G}_m$. According to Definition~\ref{GmDe}, both $G'$ and $G''$  are  isomorphic to the one-ear graph $\mathbf{P}_m$.

It follows from Conventions~\ref{conventions} that
\begin{align*}
	\RR^{V(G')}=\{\sum\limits_{i=1}^{m-1}{a_{i}\ub_{i}}
	+\sum\limits_{i=1}^{2}({b_{i}\vb_{i}+d_{i}\yb_{i}})
	+c_{1}\xb_{1}
	\mid a_{i}, b_{i},c_{1},d_{i}\in \RR \mbox{ for all } i \},
\end{align*}
where $\ub_{i}, \vb_{i}, \xb_{1}, \yb_{i}$ correspond to the vertices of $G'$ in a natural way. By the same argument as in the case of a two-ear graph, the vectors
$$
\alpha_{\ell}:=\sum\limits_{i=1}^{m-1}\ub_{i}+\xb_{1}
+2\yb_{1}+\yb_{2}
+\ell\vb_{1}+(m-\ell+1)\vb_{2},\
\ell=1,\ldots,m-1
$$
are the minimal vectors of $\mathrm{relint}(\RR_+(G'))\cap \mathcal{S}(G')$.

It is clear that $D_{G'} = \Theta_{G'}^2 y_1 v_1^{m-2} v_2^{m-1}$. The following results follow immediately from Lemmas~\ref{pd} and~\ref{canonical}.

\begin{Lemma}\label{PmNH} Assume that $G'$ is labeled as in Figure~\ref{fig3}. Then $$ \mathrm{pdim}(\mathbb{K}[G'])=m-1 \mbox{  and  }\{\Theta_{G'} v_1^{m-\ell-1}v_2^{\ell-1}\mid \ell=1,\ldots,m-1\}\subseteq\mathcal{N}_{G'}.$$
\end{Lemma}

\subsection{The complete bipartite graph $\mathbf{K}_{2,d}$}
We need to identify a subset of the top support
$\mathcal{N}_{\mathbf{K}_{2,d}}$ because the complete bipartite graph
$\mathbf{K}_{2,d}$ may also occur as an induced subgraph of $G$ for some
$d\geq 2$.

It is easy to see that $|V(\mathbf{K}_{2,d})|=d+2$.
Assume that $\mathbf{K}_{2,d}$ has two parts $\{v_1,v_2\}$ and $\{u_1,\ldots,u_d\}$.
It follows from Conventions~\ref{conventions} that
\begin{align*}
	\RR^{V(\mathbf{K}_{2,d})}
	=\{\sum\limits_{i=1}^{d}{a_{i}\ub_{i}}
	+{b_{1}\vb_{1}}+{b_{2}\vb_{2}}
	\mid a_{i}, b_{i}\in \RR \mbox{ for all } i \},
\end{align*}
where $\ub_{i},\vb_{i}$ are the unit vectors of $\RR^{V(\mathbf{K}_{2,d})}$,
each $\ub_{i}$ (resp. $\vb_{i}$) corresponds to $u_{i}$ (resp. $v_{i}$) for all $i$.

In what follows, we will construct $(d-1)$ integral vectors in $\RR^{V(\mathbf{K}_{2,d})}$ and then show that they are minimal vectors of $\mathrm{relint}(\RR_+(\mathbf{K}_{2,d}))\cap \mathcal{S}(\mathbf{K}_{2,d})$.
The construction is as follows:

For $\ell=1,\ldots,d-1$, let $$\alpha_{\ell}:=\sum\limits_{i=1}^{d}\ub_{i}+
\ell\vb_{1}+(d-\ell)\vb_{2}.$$

By \cite[Exercise~10.7.26]{V}, it follows that a vector in $\RR^{V(\mathbf{K}_{2,d})}$ of the form
\[
\sum_{i=1}^{d} a_{i}\ub_{i} + b_{1}\vb_{1} + b_{2}\vb_{2}
\]
belongs to $\RR_+(\mathbf{K}_{2,d})$ if and only if the following conditions hold:
\begin{enumerate}
\item $a_{i} \geq 0$ for all $1\leq i\leq d$;
\item $b_{i} \geq 0$ for $i=1,2$;
\item $\displaystyle\sum_{i=1}^{d} a_{i} = b_1 + b_2$.
\end{enumerate}

It is straightforward to verify that $\alpha_{\ell}$ is a minimal vector in $\mathrm{relint}(\RR_+(\mathbf{K}_{2,d}))\cap \mathcal{S}(\mathbf{K}_{2,d})$ for $\ell=1,\ldots,d-1$.
Let $H=\mathbf{K}_{2,d}$.  Then $D_H = \Theta_H^2 (v_1 v_2)^{d-2}$.  Therefore, the following result follows directly from Lemma~\ref{pd} and Lemma~\ref{canonical}.
\begin{Lemma}\label{KNH}
Let $H=\mathbf{K}_{2,d}$. Then $$\mathrm{pdim}(\mathbb{K}[H])=d-1 \mbox{  and  }\{\Theta_H v_1^{d-\ell-1}v_2^{\ell-1}\mid \ell=1,\ldots,d-1\}\subseteq\mathcal{N}_H.$$
\end{Lemma}

Thus far, we have computed the projective dimension of the edge rings for the three types of graphs under consideration. Moreover, we have determined the top support of a two-ear graph. As for the one-ear graph and the  complete bipartite graph $\mathbf{K}_{2,d}$, we have respectively obtained two particular subsets of their top supports.  We are now ready to present the first main result of this section.

\begin{Theorem} \label{third} Let $G$ be a two-ear graph of type $m$. Then, for any $0\leq i\leq m$, we have $$\beta_i(\mathbb{K}[G])=\beta_{i}(\mathbb{K}[E(G)]/\mathrm{in}_<(I_{G})).$$
 \end{Theorem}

\begin{proof} From Lemma~\ref{GmNH}, it can be seen that $\mathrm{pdim}(\mathbb{K}[G])=\mathrm{pdim}(\mathbb{K}[E(G)]/\mathrm{in}_<(I_{G}))=m$. It suffices to prove the Betti numbers of the two rings are identical.

If $i=0$, then $\beta_{i}(\mathbb{K}[G])= \beta_{i}(\mathbb{K}[E(G)]/\mathrm{in} _{<}(I_{G}))=1.$  Suppose now  that $1\leq i\leq m$. We now identify all top-Betti induced subgraphs of $G$ whose edge rings have a projective dimension of $i$. Note that $\mathrm{pdim}(\mathbb{K}[\mathbf{G}_i])=\mathrm{pdim}(\mathbb{K}[\mathbf{P}_{i+1}])=\mathrm{pdim}(\mathbb{K}[\mathbf{K}_{2,i+1}])=i$.
   Observe that $G$ has $\binom{m-1}{i-1}$ induced subgraphs that are isomorphic to $\mathbf{G}_i$,  2$\binom{m-1}{i}$ induced subgraphs that are isomorphic to $\mathbf{P}_{i+1}$, and $\binom{m-1}{i+1}$ induced subgraphs which are isomorphic to $\mathbf{K}_{2,i+1}$. Hence,   there are at least $r$ induced subgraphs $H$ of $G$ such that $\mathrm{pdim}(\mathbb{K}[H])=i$. Here, $r$ is given by
$$r=\binom{m-1}{i-1}+2\binom{m-1}{i}+\binom{m-1}{i+1}=\binom{m+1}{i+1}.$$

Let $\mathcal{N}_i(G)$ denote the disjoint union of all $\mathcal{N}_H$ such that $H$ is a top-Betti induced subgraph of $G$ satisfying $\mathrm{pdim}(\mathbb{K}[H])=i$.

Observe that if $H$ is isomorphic to $\mathbf{G}_i$, then according to Lemma~\ref{GmNH}, $\mathcal{N}_H$ comprises $i$ monomials.  If $H$ is isomorphic to either $\mathbf{P}_{i+1}$ or $\mathbf{K}_{2,i+1}$, then the corresponding lemma (Lemma~\ref{PmNH} for $\mathbf{P}_{i+1}$ and Lemma~\ref{KNH} for $\mathbf{K}_{2,i+1}$) ensures that $\mathcal{N}_H$ contains at least $i$ monomials.

 Hence,
$$|\mathcal{N}_{i}(G)|\geq ir=i\binom{m+1}{i+1}.$$
This implies
\[
i\binom{m+1}{i+1}
\leq \sum_{\alpha\in \mathcal{N}_i(G)} \beta_{i,\alpha}(\mathbb{K}[G])
\leq \beta_{i}(\mathbb{K}[G])
\leq \beta_i\bigl(\mathbb{K}[E(G)]/\mathrm{in}_{<}(I_{G})\bigr)
\leq i\binom{m+1}{i+1},
\]
where the last two inequalities follow from Lemma~\ref{Gmupper}, as desired.
\end{proof}

From the proof of Theorem~\ref{third}, we immediately obtain the following corollaries. They strengthen Lemmas~\ref{PmNH} and~\ref{KNH} by determining the full top support for
one-ear graphs and complete bipartite graphs $\mathbf{K}_{2,d}$.

\begin{Corollary}\label{PmNHeq} Let $H$ be the graph $G'$ described above. Then

$\mathrm{(1)}$ $\mathcal{N}_H=\{\Theta_H v_1^{m-\ell-1}v_2^{\ell-1}\mid \ell=1,\ldots,m-1\};$

 $\mathrm{(2)}$  the regularity of the edge ring  $\mathbb{K}[H]$ is equal to $3$.
 \end{Corollary}
\begin{Corollary}\label{PmNHeq1} Let $H$ be the graph $G''$ described above. Then

$\mathrm{(1)}$ $\mathcal{N}_H=\{\Theta_H v_1^{\ell-1}v_2^{m-\ell-1}\mid \ell=1,\ldots,m-1\};$

 $\mathrm{(2)}$  the regularity of the edge ring  $\mathbb{K}[H]$ is equal to $3$.
 \end{Corollary}

\begin{Corollary}\label{KNHeq} Let $H$ be the graph $\mathbf{K}_{2,d}$ described above. Then

$\mathrm{(1)}$ $\mathcal{N}_H=\{\Theta_H v_1^{d-\ell-1}v_2^{\ell-1}\mid \ell=1,\ldots,d-1\};$

 $\mathrm{(2)}$  the regularity of the edge ring  $\mathbb{K}[H]$ is equal to $1$.
 \end{Corollary}
We now present the second main result of this section.  As in the proof of
Theorem~\ref{third}, we use $\mathcal{N}_i(G)$ to denote the disjoint union of all
$\mathcal{N}_H$ such that $H$ is a top-Betti induced subgraph of $G$ satisfying
$\mathrm{pdim}(\mathbb{K}[H])=i$.

\begin{Corollary}\label{main2}
Let $G$ be a two-ear graph of type $m$.  Then for any monomial $\alpha$ in
$\mathbb{K}[G]$ and for any integer $i\geq 1$, we have
\[
\beta_{i,\alpha}(\mathbb{K}[G])=
\begin{cases}
1, & \text{if } \alpha\in \mathcal{N}_i(G),\\
0, & \text{otherwise.}
\end{cases}
\]

In particular, there exists a squarefree initial ideal $\mathrm{in}_{<}(I_G)$ of $I_G$
such that
\[
\beta_{i,j}(\mathrm{in}_{<}(I_G))
= \beta_{i,j}(I_G)
= \bigl|\{\alpha \in \mathcal{N}_{i+1}(G) \mid |\deg(\alpha)| = 2j\}\bigr|
\qquad\text{for all } i,j.
\]
\end{Corollary}

\begin{Example}\label{ex:main2-m3}\em
We illustrate Corollary~\ref{main2} in the case $m=3$.
Let $G=\mathbf{G}_3$ be the two-ear graph of type $3$, with vertices and edges
labeled as in Figure~1.  We focus on the multigraded Betti numbers in homological
degree $2$.

By Lemmas~\ref{GmNH}, \ref{PmNH} and~\ref{KNH} we have
\[
\mathrm{pdim}(\mathbb{K}[\mathbf{G}_2])
=\mathrm{pdim}(\mathbb{K}[\mathbf{P}_3])
=\mathrm{pdim}(\mathbb{K}[\mathbf{K}_{2,3}])
=2.
\]
Moreover, from the explicit structure of $\mathbf{G}_3$ one checks that:
\begin{itemize}
  \item $G$ has exactly two induced subgraphs isomorphic to $\mathbf{G}_2$.
        We denote these two subgraphs by $H_1$ and $H_2$.
  \item $G$ has exactly two induced subgraphs isomorphic to $\mathbf{P}_3$:
        they are the induced subgraphs $G\setminus x_2$ and $G\setminus y_2$.
        We denote them by $H_3$ and $H_4$, respectively.
\item $G$ has no induced subgraph isomorphic to $\mathbf{K}_{2,3}$.
\end{itemize}
\begin{figure}[ht]
\centering
\begin{tikzpicture}[scale=0.5,
  bigv/.style={circle,fill=black,inner sep=2pt},
  smallv/.style={circle,fill=black,inner sep=1pt}
]

\begin{scope}[shift={(0,0)}]
  \node[bigv]   (v1)  at (3,5)   {};
  \node[bigv]   (v2)  at (9,5)   {};
  \node[smallv] (x1)  at (3,1)   {};
  \node[smallv] (u2)  at (6,3)   {};
  \node[smallv] (y1)  at (9,1)   {};
  \node[smallv] (u1)  at (6,7)   {};
  \node[smallv] (x2)  at (1.3,3) {};
  \node[smallv] (y2)  at (10.7,3){};

  \draw[black,thin] (x1)--(v1)--(u2)--(v2)--(y1)--(x1);  
  \draw[black,thin] (x1)--(x2)--(v1)--(x1);             
  \draw[black,thin] (y1)--(y2)--(v2)--(y1);             
  \draw[black,thin,dashed] (v1)--(u1);
  \draw[black,thin,dashed] (u1)--(v2);

  \node[above]  at (v1) {$v_1$};
  \node[above]  at (v2) {$v_2$};
  \node[below]  at (x1) {$x_1$};
  \node[below]  at (y1) {$y_1$};
  \node[right]  at (u2) {$u_2$};
  \node[above]  at (u1) {$u_1$};
  \node[left]   at (x2) {$x_2$};
  \node[right]  at (y2) {$y_2$};

  \node[below=10pt] at (6,0) {$H_1$};
\end{scope}

\begin{scope}[shift={(15,0)}]
  \node[bigv]   (v1)  at (3,5)   {};
  \node[bigv]   (v2)  at (9,5)   {};
  \node[smallv] (x1)  at (3,1)   {};
  \node[smallv] (u2)  at (6,3)   {};
  \node[smallv] (y1)  at (9,1)   {};
  \node[smallv] (u1)  at (6,7)   {};
  \node[smallv] (x2)  at (1.3,3) {};
  \node[smallv] (y2)  at (10.7,3){};

  \draw[black,thin] (x1)--(v1)--(u1)--(v2)--(y1)--(x1);  
  \draw[black,thin] (x1)--(x2)--(v1)--(x1);             
  \draw[black,thin] (y1)--(y2)--(v2)--(y1);             
  \draw[black,thin,dashed] (v1)--(u2);
  \draw[black,thin,dashed] (u2)--(v2);

  \node[above]  at (v1) {$v_1$};
  \node[above]  at (v2) {$v_2$};
  \node[below]  at (x1) {$x_1$};
  \node[below]  at (y1) {$y_1$};
  \node[right]  at (u2) {$u_2$};
  \node[above]  at (u1) {$u_1$};
  \node[left]   at (x2) {$x_2$};
  \node[right]  at (y2) {$y_2$};

  \node[below=10pt] at (6,0) {$H_2$};
\end{scope}

\begin{scope}[shift={(0,-12)}]
  \node[bigv]   (v1)  at (3,5)   {};
  \node[bigv]   (v2)  at (9,5)   {};
  \node[smallv] (x1)  at (3,1)   {};
  \node[smallv] (u2)  at (6,3)   {};
  \node[smallv] (y1)  at (9,1)   {};
  \node[smallv] (u1)  at (6,7)   {};
  \node[smallv] (x2)  at (1.3,3) {};
  \node[smallv] (y2)  at (10.7,3){};

  \draw[black,thin] (x1)--(v1)--(u2)--(v2)--(y1)--(x1);  
  \draw[black,thin] (y1)--(y2)--(v2)--(y1);             
  \draw[black,thin] (v1)--(u1)--(v2);                   
  \draw[black,thin,dashed] (x1)--(x2);
  \draw[black,thin,dashed] (x2)--(v1);

  \node[above]  at (v1) {$v_1$};
  \node[above]  at (v2) {$v_2$};
  \node[below]  at (x1) {$x_1$};
  \node[below]  at (y1) {$y_1$};
  \node[right]  at (u2) {$u_2$};
  \node[above]  at (u1) {$u_1$};
  \node[left]   at (x2) {$x_2$};
  \node[right]  at (y2) {$y_2$};

  \node[below=10pt] at (6,0) {$H_3$};
\end{scope}

\begin{scope}[shift={(15,-12)}]
  \node[bigv]   (v1)  at (3,5)   {};
  \node[bigv]   (v2)  at (9,5)   {};
  \node[smallv] (x1)  at (3,1)   {};
  \node[smallv] (u2)  at (6,3)   {};
  \node[smallv] (y1)  at (9,1)   {};
  \node[smallv] (u1)  at (6,7)   {};
  \node[smallv] (x2)  at (1.3,3) {};
  \node[smallv] (y2)  at (10.7,3){};

  \draw[black,thin] (x1)--(v1)--(u2)--(v2)--(y1)--(x1);  
  \draw[black,thin] (x1)--(x2)--(v1)--(x1);             
  \draw[black,thin] (v1)--(u1)--(v2);                   
  \draw[black,thin,dashed] (y1)--(y2);
  \draw[black,thin,dashed] (y2)--(v2);

  \node[above]  at (v1) {$v_1$};
  \node[above]  at (v2) {$v_2$};
  \node[below]  at (x1) {$x_1$};
  \node[below]  at (y1) {$y_1$};
  \node[right]  at (u2) {$u_2$};
  \node[above]  at (u1) {$u_1$};
  \node[left]   at (x2) {$x_2$};
  \node[right]  at (y2) {$y_2$};

  \node[below=10pt] at (6,0) {$H_4$};
\end{scope}

\end{tikzpicture}
\caption{The graphs $H_1,H_2,H_3,H_4$.}
\label{fig:H1234}
\end{figure}
By Proposition~\ref{prop:two-ear-top-betti-types}, the top-Betti induced
subgraphs $H$ of $G$ with $\mathrm{pdim}(\mathbb{K}[H])=2$ are therefore
precisely the four graphs $H_1,H_2,H_3,H_4$.  Consequently,
\[
\begin{aligned}
\mathcal{N}_2(G)
&= \mathcal{N}_{H_1} \sqcup \mathcal{N}_{H_2}
   \sqcup \mathcal{N}_{H_3} \sqcup \mathcal{N}_{H_4} \\[4pt]
&= \bigl\{
   \Theta_{H_1} x_1y_1v_1,\;
   \Theta_{H_1} x_1y_1v_2,\;
   \Theta_{H_2} x_1y_1v_1,\;
   \Theta_{H_2} x_1y_1v_2,
\\
&\qquad\quad
   \Theta_{H_3} v_1,\;
   \Theta_{H_3} v_2,\;
   \Theta_{H_4} v_1,\;
   \Theta_{H_4} v_2
   \bigr\} \\[4pt]
&= \bigl\{
   x_1^2 y_1^2 x_2 y_2 v_1^2 v_2 u_2,\;
   x_1^2 y_1^2 x_2 y_2 v_1 v_2^2 u_2,\;
   x_1^2 y_1^2 x_2 y_2 v_1^2 v_2 u_1,\;
   x_1^2 y_1^2 x_2 y_2 v_1 v_2^2 u_1,
\\
&\qquad\quad
   x_1 y_1 y_2 v_1^2 v_2 u_1u_2,\;
   x_1 y_1 y_2 v_1 v_2^2 u_1u_2,\;
   x_1 y_1 x_2 v_1 v_2^2 u_1u_2,\;
   x_1 y_1 x_2 v_1^2 v_2 u_1u_2
   \bigr\}.
\end{aligned}
\]

Applying Corollary~\ref{main2} with $i=2$ (and $m=3$), we obtain the multigraded and graded Betti numbers of $\mathbb{K}[G]$ in homological degree 2:
\[
\beta_{2,\alpha}(\mathbb{K}[G])=
\begin{cases}
1, & \text{if } \alpha \in \mathcal{N}_2(G),\\
0, & \text{otherwise,}
\end{cases}
\qquad
\beta_{2,j}(\mathbb{K}[G])=
\begin{cases}
4, & \text{if } j=4 \text{ or } j=5,\\
0, & \text{otherwise.}
\end{cases}
\]

\end{Example}

\section{Compact graphs revisited}
The primary objective of this section is to determine the elements of the top support $\mathcal{N}_H$ for every compact graph $H$, thereby facilitating the application of the induced-subgraph approach for computing the Betti numbers of their edge rings in the subsequent section.  To accomplish this, we revisit the precise definition and classification of compact graphs as presented in \cite{WL}.

\begin{Definition}\label{3.1}\em Let $G$ be a connected simple graph where every vertex has degree at least $2$. We call $G$ a {\it compact} graph if it does not contain any even cycles and satisfies the odd cycle condition.
\end{Definition}

A vertex of a compact graph is called {\it big} if its degree is greater
than $2$. According to \cite[Lemma~3.10]{WL}, there are at most three
big vertices in a compact graph. We say that a compact graph is of
type~$i$ if it has exactly $i$ big vertices, where $0\leq i\leq3$.
\begin{Conventions}\em
	A compact graph of type~0 is an odd cycle. In this case the edge ring is a
	polynomial ring, so its Betti numbers are trivial. Therefore, in the sequel,
	the term {\it compact graph} excludes odd cycles unless explicitly stated
	otherwise.
\end{Conventions}

We now present a complete classification of compact graphs, originally established in \cite[Theorem~3.12]{WL}.
A compact graph of type 1 is obtained from at least two odd cycles by identifying one vertex from each cycle to a single common vertex.

A compact graph of type 2 is obtained from two vertex-disjoint graphs, each formed from one or more odd cycles by identifying one vertex from each cycle to a distinguished vertex, by connecting the two distinguished vertices either by an edge or by an edge together with a path of even length.

A compact graph of type 3 is obtained from three pairwise vertex-disjoint graphs, each formed from one or more odd cycles by identifying one vertex from each cycle to a distinguished vertex, by joining every pair of the three distinguished vertices by an edge.

\begin{figure}[ht]
\centering
\begin{tikzpicture}[
  scale=1.0,
  bigv/.style={circle, fill=black, inner sep=1.5pt},
  smallv/.style={circle, fill=black, inner sep=0.9pt}
]

\begin{scope}[shift={(-0.5,0)}]
  \node[bigv] (c) at (0,0) {};

  \node[draw=none, fill=none] at (0,0.28) {$u$};

  \node[smallv] (a1) at (1,0.8) {};
  \node[smallv] (b1) at (1.4,0) {};
  \draw (c)--(a1)--(b1)--(c);

  \node[smallv] (a3) at (-1,0.8) {};
  \node[smallv] (b3) at (-1.4,0) {};
  \draw (c)--(a3)--(b3)--(c);

  \node[smallv] (p2) at (-0.761,-0.553) {};
  \node[smallv] (p3) at (-0.47,-1.447) {};
  \node[smallv] (p4) at (0.47,-1.447) {};
  \node[smallv] (p5) at (0.761,-0.553) {};
  \draw (c)--(p5)--(p4)--(p3)--(p2)--(c);

  \node[draw=none, fill=none] at (0,-2.4) {\small $\mathbf{A}_{(1,1,2)}$};
\end{scope}

\begin{scope}[shift={(2.9,0)}]  
  \node[bigv] (L)  at (0,0) {};
  \node[smallv] (L1) at (-0.9,0.8) {};
  \node[smallv] (L2) at (-0.9,-0.8) {};
  \draw (L)--(L1)--(L2)--(L);

  \node[draw=none, fill=none] at (0,-0.35) {$u$};

  \node[bigv] (R)  at (2.6,0) {};
  \node[smallv] (R1) at (3.5,0.8) {};
  \node[smallv] (R2) at (3.5,-0.8) {};
  \draw (R)--(R1)--(R2)--(R);

  \node[draw=none, fill=none] at (2.6,-0.35) {$v$};

  \draw (L)--(R);             

  \node[smallv] (M1) at (0.65,0.9) {};
  \node[smallv] (M2) at (1.30,0.9) {};
  \node[smallv] (M3) at (1.95,0.9) {};
  \draw (L)--(M1)--(M2)--(M3)--(R); 

  \node[draw=none, fill=none] at (1.3,-2.4) {\small $\mathbf{B}^{4}_{(1):(1)}$};
\end{scope}

\begin{scope}[shift={(8.1,-0.6)}]  
  \node[bigv] (B1) at (0,0) {};
  \node[bigv] (B2) at (2.4,0) {};
  \node[bigv] (B3) at (1.2,2.0) {};

  \node[draw=none, fill=none] at (0,-0.35) {$u$};
  \node[draw=none, fill=none] at (2.4,-0.35) {$v$};
  \node[draw=none, fill=none] at (1.55,2.0) {$w$};

  \draw (B1)--(B2)--(B3)--(B1);


  \node[smallv] (B1a) at (-0.9,0.8) {};
  \node[smallv] (B1b) at (-0.9,-0.8) {};
  \draw (B1)--(B1a)--(B1b)--(B1);

  \node[smallv] (B2a) at (3.3,0.8) {};
  \node[smallv] (B2b) at (3.3,-0.8) {};
  \draw (B2)--(B2a)--(B2b)--(B2);

  \node[smallv] (B3a) at (0.6,3.0) {};
  \node[smallv] (B3b) at (1.8,3.0) {};
  \draw (B3)--(B3a)--(B3b)--(B3);

  \node[draw=none, fill=none] at (1.2,-2.4) {\small $\mathbf{C}_{(1):(1):(1)}$};
\end{scope}

\end{tikzpicture}
\caption{Examples of compact graphs of types 1, 2 and 3.}
\end{figure}

Suppose $\underline{p}=(p_1,\ldots,p_m)$, $\underline{q}=(q_1,\ldots,q_n)$, and $\underline{r}=(r_1,\ldots,r_k)$ are positive integral vectors of dimensions $m,n$, and $k$, respectively.
We denote by $\mathbf{A}_{\underline{p}}$ (or
$\mathbf{A}_{(p_1,\ldots,p_m)}$) the graph obtained from $m$ odd cycles
of lengths $2p_1+1,\ldots,2p_m+1$, respectively, by identifying one
vertex from each cycle to a single distinguished vertex. Thus,
$\mathbf{A}_{\underline{p}}$ is a compact graph of type~0 when $m=1$
and of type~1 when $m\geq2$; in the latter case, its distinguished
vertex is its unique big vertex.

By $\mathbf{B}_{\underline{p}:\underline{q}}^0$, we mean the compact graph of type~2 obtained from vertex-disjoint copies of $\mathbf{A}_{\underline{p}}$ and $\mathbf{A}_{\underline{q}}$ by joining their distinguished vertices by an edge. Furthermore, if $s\geq2$ is even, then $\mathbf{B}_{\underline{p}:\underline{q}}^s$ denotes the graph obtained from $\mathbf{B}_{\underline{p}:\underline{q}}^0$ by appending a path of length $s$ between the two distinguished vertices.

We denote by $\mathbf{C}_{\underline{p}:\underline{q}:\underline{r}}$ the compact graph of type~3 obtained from pairwise vertex-disjoint copies of $\mathbf{A}_{\underline{p}}$, $\mathbf{A}_{\underline{q}}$, and $\mathbf{A}_{\underline{r}}$ by joining every pair of their distinguished vertices by an edge.

We recall the results from \cite{WL} concerning the minimal vectors of
$\operatorname{relint}(\RR_+(G))\cap\mathcal{S}(G)$ for compact graphs,
and subsequently determine the top supports of the three types of compact
graphs, respectively. Throughout this section, we use the notation for $\RR^{V(G)}$ introduced
in Conventions~\ref{conventions}.

\subsection{Type one}

Given an integer $m \geq 2$ and a positive integer tuple $\underline{p} =(p_{1}, \ldots, p_{m})$, we denote $\mathbf{A}_{\underline{p}}$ simply by $A$.  We may assume that $A$ has vertex set
\[
V(A)=\{u\} \sqcup \{u_{i,j} \mid 1 \leq i \leq m, 1 \leq j \leq 2p_{i}\}
\]
and edge set
\begin{align*}
E(A)=\{& \{u_{i,j},u_{i,j+1}\} \mid 1 \leq i \leq m, 1 \leq j \leq 2p_{i}-1 \}\\
&\sqcup \{ \{u,u_{i,1}\}, \{u,u_{i,2p_{i}}\} \mid  1 \leq i \leq m \}.
\end{align*}
 Recall again that an integral vector in $\mathrm{relint}(\RR_+(G))\cap \mathcal{S}(G)$ is   minimal if it cannot be written as the sum of a vector in $\mathrm{relint}(\RR_+(G))\cap \mathcal{S}(G)$ and a nonzero vector of  $\RR_+(G)\cap \mathcal{S}(G)$.
We  showed in \cite[Section 6.2]{WL} that the following vectors
 $$\alpha_{\ell}:=\sum\limits_{i=1}^{m}\sum\limits_{j=1}^{2p_{i}}\ub_{i,j}+2\ell \ub,\  \ell=1,\ldots,m-1$$
 are precisely the minimal vectors of $\mathrm{relint}(\RR_+(A))\cap \mathcal{S}(A)$.

It is clear that $D_A=\Theta_A^2 u^{2m-2}$.  Therefore, the following result follows directly from Lemma~\ref{canonical}.

\begin{Lemma} \label{5.2} Let $H$ be a compact graph  of type $1$ labeled as previously. Then $$\mathcal{N}_H=\{\Theta_H u^{2(m-\ell)-1}\mid \ell=1,\ldots,m-1\}. $$

\end{Lemma}

\subsection{Type two}

Given positive integers $m,n$, an even integer $s\geq 0$, and positive
integral tuples $\underline{p}=(p_1,\ldots,p_m)$ and
$\underline{q}=(q_1,\ldots,q_n)$, we denote
$\mathbf{B}_{\underline{p}:\underline{q}}^s$ simply by $B^s$.

Set
\[
W_s=
\begin{cases}
	\varnothing, & s=0,\\
	\{w_1,\ldots,w_{s-1}\}, & s\geq2,
\end{cases}
\]
and let
\[
P_s=
\begin{cases}
	\bigl\{\{u,v\}\bigr\}, & s=0,\\[2mm]
	\bigl\{\{u,v\},\{u,w_1\},\{v,w_{s-1}\}\bigr\}
	\cup
	\bigl\{\{w_i,w_{i+1}\}\mid 1\leq i\leq s-2\bigr\},
	& s\geq2.
\end{cases}
\]
Thus, when $s=0$, the two distinguished vertices $u$ and $v$ are
joined only by the edge $\{u,v\}$, whereas when $s\geq2$, they are
joined by the edge $\{u,v\}$ together with a path of length $s$.

The vertex set of $B^s$ is
\[
V(B^s)
=
\{u,v\}
\sqcup W_s
\sqcup
\{u_{i,j}\mid 1\leq i\leq m,\ 1\leq j\leq2p_i\}
\sqcup
\{v_{i,j}\mid 1\leq i\leq n,\ 1\leq j\leq2q_i\},
\]
and its edge set is
\begin{align*}
	E(B^s)
	={}&
	\bigl\{
	\{u_{i,j},u_{i,j+1}\}
	\mid
	1\leq i\leq m,\ 1\leq j\leq2p_i-1
	\bigr\}
	\\
	&\sqcup
	\bigl\{
	\{u,u_{i,1}\},\{u,u_{i,2p_i}\}
	\mid
	1\leq i\leq m
	\bigr\}
	\\
	&\sqcup P_s
	\\
	&\sqcup
	\bigl\{
	\{v_{i,j},v_{i,j+1}\}
	\mid
	1\leq i\leq n,\ 1\leq j\leq2q_i-1
	\bigr\}
	\\
	&\sqcup
	\bigl\{
	\{v,v_{i,1}\},\{v,v_{i,2q_i}\}
	\mid
	1\leq i\leq n
	\bigr\}.
\end{align*}

We showed in \cite[Section~6.3]{WL} that the following vectors
 $$\alpha_{\ell}:=\sum\limits_{i=1}^{m}\sum\limits_{j=1}^{2p_{i}}\ub_{i,j}+\sum\limits_{i=1}^{n}\sum\limits_{j=1}^{2q_{i}}\vb_{i,j}+\vb+(2\ell+1)\ub, \quad \ell=0,\ldots,m-1$$
and $$\beta_{\ell}:=\sum\limits_{i=1}^{m}\sum\limits_{j=1}^{2p_{i}}\ub_{i,j}+\sum\limits_{i=1}^{n}\sum\limits_{j=1}^{2q_{i}}\vb_{i,j}+\ub+(2\ell+1)\vb, \quad \ell=1,\ldots,n-1$$
are precisely the minimal vectors of $\mathrm{relint}(\RR_+(B^0))\cap \mathcal{S}(B^0)$.

Similarly, for any even integer $s\geq2$, it was shown in
\cite[Section~6.3]{WL} that the vectors
$$
\alpha_{\ell}:=
\sum_{i=1}^{m}\sum_{j=1}^{2p_{i}}\ub_{i,j}
+\sum_{i=1}^{n}\sum_{j=1}^{2q_{i}}\vb_{i,j}
+\sum_{i=1}^{s-1}\wb_{i}
+\vb+2\ell\ub,\qquad \ell=1,\ldots,m,
$$
and
$$
\beta_{\ell}:=
\sum_{i=1}^{m}\sum_{j=1}^{2p_{i}}\ub_{i,j}
+\sum_{i=1}^{n}\sum_{j=1}^{2q_{i}}\vb_{i,j}
+\sum_{i=1}^{s-1}\wb_{i}
+\ub+2\ell\vb,\qquad \ell=1,\ldots,n,
$$
are precisely the minimal vectors of
$\mathrm{relint}(\RR_+(B^s))\cap \mathcal{S}(B^s)$.

\begin{Lemma} \label{5.3}
For each compact graph $H$ of type 2,
\begin{enumerate}
  \item[$(1)$] If $H=B^0$, then
  $$\mathcal{N}_H=\{\Theta_H u^{2(m-\ell)-1}v^{2n-1}\mid \ell=0,\ldots,m-1\}\sqcup  \{\Theta_H u^{2m-1}v^{2(n-\ell)-1}\mid \ell=1,\ldots,n-1\}.$$
  \item[$(2)$] If $H=B^s$ for some even integer $s\geq 2$, then
  $$\mathcal{N}_H=\{\Theta_H u^{2(m-\ell)+1}v^{2n}\mid \ell=1,\ldots,m\}\sqcup  \{\Theta_H u^{2m}v^{2(n-\ell)+1}\mid \ell=1,\ldots,n\}.$$
\end{enumerate}

\begin{proof}
If $H=B^0$, then $D_H = \Theta_{H}^2 u^{2m-1}v^{2n-1}$.
Combining the above description of the minimal vectors of
$\mathrm{relint}(\RR_+(B^0))\cap \mathcal{S}(B^0)$ with
Lemma~\ref{canonical} yields the formula for $\mathcal{N}_H$ in \emph{(1)}.

Similarly, if $H=B^s$ for some even integer $s\ge 2$, then
$D_H = \Theta_{H}^2 u^{2m}v^{2n}$.
Combining the above description of the minimal vectors of
$\mathrm{relint}(\RR_+(B^s))\cap \mathcal{S}(B^s)$ with
Lemma~\ref{canonical} yields the formula for $\mathcal{N}_H$ in \emph{(2)}.
\end{proof}

\end{Lemma}

In Lemma \ref{5.3}, the notation $u,v$ used in the expression of $\mathcal{N}_H$ denotes the two big vertices of the graph $H$. Nevertheless, when $H$ is regarded as an induced subgraph of a compact graph of type~3, these two big vertices may instead be denoted by $w,u$ or $w,v$. In such cases, one has to rename the vertices so that the pair of big vertices is written as $u,v$ in accordance with the notation used in Lemma~\ref{5.3}.

\subsection{Type three}
Given positive integers $m,n,k$, as well as the positive integral  tuples $\underline{p}=(p_{1}, \ldots, p_{m})$, $\underline{q}=(q_{1}, \ldots, q_{n})$ and  $\underline{r}=(r_{1}, \ldots, r_{k})$,
we denote $\mathbf{C}_{\underline{p}:\underline{q}:\underline{r}}$ simply by $C$. By definition, we may assume that $C$ has  vertex set
\begin{align*}
V(C)=&\{ u,v,w\} \sqcup \{u_{i,j}\mid 1 \leq i \leq m, 1 \leq j \leq 2p_{i}\}\\
& \sqcup \{v_{i,j}\mid 1 \leq i \leq n, 1 \leq j \leq 2q_{i}\} \sqcup \{w_{i,j}\mid 1 \leq i \leq k, 1 \leq j \leq 2r_{i}\}
\end{align*}

and edge set
\begin{align*}
E(C)=&\{ \{u_{i,j},u_{i,j+1}\} \mid 1 \leq i \leq m, 1 \leq j \leq 2p_{i}-1 \}\\
&\sqcup \{ \{u,u_{i,1}\}, \{u,u_{i,2p_{i}}\} \mid  1 \leq i \leq m \} \\
&\sqcup \{ \{v_{i,j},v_{i,j+1}\} \mid 1 \leq i \leq n, 1 \leq j \leq 2q_{i}-1 \}\\
&\sqcup \{ \{v,v_{i,1}\}, \{v,v_{i,2q_{i}}\} \mid  1 \leq i \leq n \}\\
&\sqcup \{ \{w_{i,j},w_{i,j+1}\} \mid 1 \leq i \leq k, 1 \leq j \leq 2r_{i}-1 \}\\
&\sqcup \{ \{w,w_{i,1}\}, \{w,w_{i,2r_{i}}\} \mid  1 \leq i \leq k \}\\
&\sqcup \{\{u,v\}, \{v, w\},\{w,u\}\} .
\end{align*}

 By \cite[Section 6.1]{WL}, we have the following result.

\begin{Lemma} \label{3.4}
The minimal vectors of $\mathrm{relint}(\RR_+(C))\cap \mathcal{S}(C)$ are precisely the following $m+n+k$ vectors:
$\alpha_{\ell}$ for $\ell=1,\ldots,m$, $\beta_{\ell}$ for $\ell=1,\ldots,n$, and $\gamma_{\ell}$ for $\ell=1,\ldots,k$.

For $\ell=1,\ldots,m$, set
$$
\alpha_{\ell}:=
\sum_{i=1}^{m}\sum_{j=1}^{2p_{i}}\ub_{i,j}
+\sum_{i=1}^{n}\sum_{j=1}^{2q_{i}}\vb_{i,j}
+\sum_{i=1}^{k}\sum_{j=1}^{2r_{i}}\wb_{i,j}
+\vb+\wb+2\ell\ub.
$$

For $\ell=1,\ldots,n$, set
$$
\beta_{\ell}:=
\sum_{i=1}^{m}\sum_{j=1}^{2p_{i}}\ub_{i,j}
+\sum_{i=1}^{n}\sum_{j=1}^{2q_{i}}\vb_{i,j}
+\sum_{i=1}^{k}\sum_{j=1}^{2r_{i}}\wb_{i,j}
+\wb+\ub+2\ell\vb.
$$

For $\ell=1,\ldots,k$, set
$$
\gamma_{\ell}:=
\sum_{i=1}^{m}\sum_{j=1}^{2p_{i}}\ub_{i,j}
+\sum_{i=1}^{n}\sum_{j=1}^{2q_{i}}\vb_{i,j}
+\sum_{i=1}^{k}\sum_{j=1}^{2r_{i}}\wb_{i,j}
+\ub+\vb+2\ell\wb.
$$

In other words, the canonical module $\omega_{\mathbb{K}[C]}$ of $\mathbb{K}[C]$ is a monomial ideal in $\mathbb{K}[C]$ generated by
$$
\Theta_C u, \ \Theta_C u^3, \ \ldots, \ \Theta_C u^{2m-1}, \quad
\Theta_C v, \ \Theta_C v^3, \ \ldots, \ \Theta_C v^{2n-1},
$$
and
$$
\Theta_C w, \ \Theta_C w^3, \ \ldots, \ \Theta_C w^{2k-1}.
$$
\end{Lemma}

The computation of $D_C$ is straightforward, resulting in $D_C=\Theta_C^2 u^{2m}v^{2n}w^{2k}$. Applying Lemma~\ref{canonical} to this equality, we obtain the following result immediately.
\begin{Lemma}\label{5.5}

If $H=\mathbf{C}_{\underline{p}:\underline{q}:\underline{r}}$, then \begin{align*}\mathcal{N}_H=&\{\Theta_H u^{2(m-\ell)+1}v^{2n}w^{2k}\mid \ell=1,\ldots,m\}\sqcup  \{\Theta_H u^{2m}v^{2(n-\ell)+1}w^{2k}\mid  \ell=1,\ldots,n\}\\ &\sqcup \{\Theta_H u^{2m}v^{2n}w^{2(k-\ell)+1}\mid \ell=1,\ldots,k\}.
\end{align*}
\end{Lemma}

So far, we have computed the top supports of all compact graphs. From the
above formulas, we conclude that, for each compact graph $H$, its top support
$\mathcal{N}_H$ contains precisely $\mathfrak{t}(H)-1$ elements, where
$\mathfrak{t}(H)$ denotes the number of minimal cycles in $H$.
Moreover, for any compact graph $H$, we have that $\Theta_H$ divides $\alpha$
for every $\alpha\in \mathcal{N}_H$. This shows that every compact graph is a
top-Betti graph, which, by Proposition~\ref{second-Betti}, in turn implies that
every compact graph is also a second-Betti graph. Alternatively, a proof of this
conclusion can also be obtained by applying Proposition~\ref{Betti}.

We now describe the top-Betti induced subgraphs of a compact graph.

\begin{Proposition}\label{compact-top-Betti}
Let $G$ be a compact graph. Then the induced subgraphs of $G$ that are top-Betti
graphs are precisely the compact induced subgraphs of $G$.
\end{Proposition}

\begin{proof}
Let $H$ be an induced subgraph of $G$. If $H$ is compact, then $H$ is top-Betti by the
preceding discussion. Thus it remains to show that every non-compact induced subgraph
of $G$ is not top-Betti.

By the classification of compact graphs, every induced subgraph \( H \) of a compact graph \( G \) is either compact, an odd cycle, or one with a vertex of degree at most one. By Corollary~\ref{top-Betti at least two}, \( H \) is not a top-Betti graph if it contains a vertex of degree at most one. Similarly, Proposition~\ref{Betti} implies that an odd cycle cannot be a top-Betti graph.

Therefore the only induced subgraphs of $G$ that are top-Betti graphs are the compact
ones, as required.
\end{proof}

\section{Betti numbers of compact graphs}
In this section, we employ the induced-subgraph approach to determine
the multigraded Betti numbers of the edge rings of compact graphs.

\begin{Lemma} \label{compactupper} Let $G$ be a compact graph. Then

$\mathrm{(1)}$ The projective dimension of $\mathbb{K}[G]$ is equal to  $\mathfrak{t}(G) - 1$;

$\mathrm{(2)}$ There is a monomial order $<$ such  that   $\beta_{i}(\mathrm{in} _{<}(I_{G}))=(i+1)\binom{\mathfrak{t}(G)}{i+2}$ for all $i\geq 0$.  Thus, for all $i\geq 1$, one has $$\beta_i(\mathbb{K}[G])\leq \beta_i(\mathbb{K}[E(G)]/\mathrm{in} _{<}(I_{G}))=i\binom{\mathfrak{t}(G)}{i+1}.$$
\begin{proof} The conclusion (1) follows from \cite [Corollary 5.2]{WL}, and the conclusion (2) follows from \cite [Theorem 5.1]{WL} as well as Lemma~\ref {total}.
\end{proof}
\end{Lemma}

In the subsequent discussions, the monomial order denoted by \( < \) always refers to the order introduced in Lemma~\ref{compactupper}, which was first defined in \cite[Section 4]{WL}. For our purposes, an explicit description of this order is unnecessary, so we omit it.

\subsection{Type 1 and 2}
\begin{Proposition}\label{first} Let $G$ be a compact graph of either type 1 or type 2. Then, for any $0\leq i\leq \mathfrak{t}(G)-1$, we have $$\beta_i(\mathbb{K}[G])=\beta_{i}(\mathbb{K}[E(G)]/\mathrm{in}_<(I_G)) .$$
 \end{Proposition}

\begin{proof} If $i=0$, then $\beta_{i}(\mathbb{K}[G])= \beta_{i}(\mathbb{K}[E(G)]/\mathrm{in} _{<}(I_{G}))=1.$  Fix $1\leq i\leq \mathfrak{t}(G)-1$. Since there are $\mathfrak{t}(G)$ minimal cycles, there are $\binom{\mathfrak{t}(G)}{i+1}$ compact induced subgraphs of $G$ that contain  exactly $i+1$ minimal cycles of $G$. These induced subgraphs are of either type one or type two, and we denote them by $H_1,\ldots, H_r$. Here, $r=\binom{\mathfrak{t}(G)}{i+1}$. Note that $\beta_{i,\alpha}(\mathbb{K}[H_\ell])=1$ for any $\alpha\in \mathcal{N}_{H_\ell}$ by Lemmas~\ref{5.2} and~\ref{5.3}. It follows that $\beta_{i,\alpha}(\mathbb{K}[G])=1$ for any $\alpha\in \bigcup_{\ell=1}^r\mathcal{N}_{H_{\ell}}$ in view of Lemma~\ref{start}. Put $\mathcal{N}_i(G)=\bigcup_{\ell=1}^r\mathcal{N}_{H_{\ell}}$.

Since $\mathrm{supp}(\alpha)=V(H_k)$ for any $\alpha\in \mathcal{N}_{H_k}$, the sets $\mathcal{N}_{H_1},\ldots, \mathcal{N}_{H_r}$ are pairwise disjoint. Hence,
$$|\mathcal{N}_i(G)|=\sum_{\ell=1}^r|\mathcal{N}_{H_\ell}|=ir=i\binom{\mathfrak{t}(G)}{i+1}.$$ This implies $$i\binom{\mathfrak{t}(G)}{i+1}=\sum_{\alpha\in \mathcal{N}_i(G)}\beta_{i,\alpha}(\mathbb{K}[G])\leq\beta_{i}(\mathbb{K}[G])\leq \beta_i(\mathbb{K}[E(G)]/\mathrm{in} _{<}(I_{G}))=i\binom{\mathfrak{t}(G)}{i+1},$$
completing the proof.
\end{proof}

From the proof of Proposition~\ref{first}, it is evident that  $$\beta_{i}(\mathbb{K}[G])=\sum_{\alpha\in \mathcal{N}_i(G)}\beta_{i,\alpha}(\mathbb{K}[G]).$$ In particular,  if $\beta_{i,\alpha}(\mathbb{K}[G])\neq 0$, then $\alpha\in \mathcal{N}_i(G)$. Therefore, we can derive the multigraded  and graded Betti numbers of $\mathbb{K}[G]$ as follows.

\begin{Corollary} \label{compactBetti12}
	Let $G$ be a compact graph of either type 1 or type 2.
	For every integer $1\leq i\leq \mathfrak{t}(G)-1$ and every monomial
	$\alpha$ in $\mathbb{K}[G]$, we have $$\beta_{i,\alpha}(\mathbb{K}[G])=\left\{
                                          \begin{array}{ll}
                                            1, & \hbox{$\alpha\in \mathcal{N}_i(G)$;} \\
                                            0, & \hbox{otherwise.}
                                          \end{array}
                                        \right.$$
 In particular, $\beta_{i,j}(\mathbb{K}[G])$ is equal to the cardinality
 of the set $\{\alpha\in \mathcal{N}_i(G)\mid |\deg(\alpha)|=2j\}$.
 \end{Corollary}

 Corollary~\ref{compactBetti12} implies that if $G$ is a compact graph of type 1 or type 2 then $G$ is both a top-Betti graph and a second-Betti graph, but not a proper second-Betti graph.

We now turn to compact graphs of type 3. In the case when $G$ is of type 1 or type 2,  there are   exactly $\binom{\mathfrak{t}(G)}{i+1}$  distinct compact induced subgraphs of $G$ which contain exactly $i+1$ minimal cycles for all $1\leq i\leq \mathfrak{t}(G)-1$. This is the reason why the proof of Proposition~\ref{first} holds. But this is not the case when considering  compact graphs of type 3. Consider the following example.

\begin{Example} \label{example1}\em
Assume that $G=\mathbf{C}_{(1):(1):(1,1)}$, as shown in Figure~\ref{fig}. Then $G$ contains 5  induced cycles. But there are only 4 distinct compact induced subgraphs of $G$  each containing 4 induced cycles. This is because the union of the 4 cycles around the triangle formed by the three big vertices $\{u,v,w\}$ does not form an induced subgraph of $G$. Thus, unlike the type~1 and type~2 cases, the top-Betti contributions
do not attain the upper bound in homological degree~$3$. This suggests
that second-top contributions may be needed.

\begin{figure}[ht]
\centering
\begin{tikzpicture}

\draw[black, thin] (5,2) -- (3,2.2) -- (3.3,3.2)-- cycle;
\draw[black, thin] (5,2) -- (6.1,4) -- (7.2,2)-- cycle;
\draw[black, thin] (7.2,2) -- (9.2,2.2) -- (8.9,3.2)-- cycle;

\draw[black, thin] (6.1,4) -- (4,5.5) -- (3.5,4.5)-- cycle;
\draw[black, thin] (6.1,4) -- (8.2,5.5) -- (8.7,4.5)-- cycle;

\filldraw [black] (5,2) circle (1.5pt);
\filldraw [black] (7.2,2) circle (1.5pt);
\filldraw [black] (6.1,4) circle (1.5pt);

\draw (8.2,5.4) node[anchor=south]{$x_1$};
\draw (8.9,4.2) node[anchor=south]{$x_2$};

\draw (4.1,5.4) node[anchor=south]{$w_1$};
\draw (3.4,4.0) node[anchor=south]{$w_2$};
\draw (3.4,3.1) node[anchor=south]{$u_2$};
\draw (2.9,1.7) node[anchor=south]{$u_1$};

\draw (9.3,1.7) node[anchor=south]{$v_1$};
\draw (8.9,3.2) node[anchor=south]{$v_2$};

\filldraw [black] (8.2,5.5) circle (0.9pt);
\filldraw [black] (8.7,4.5) circle (0.9pt);

\filldraw [black] (4,5.5) circle (0.9pt);
\filldraw [black] (3.5,4.5) circle (0.9pt);

\filldraw [black] (3,2.2) circle (0.9pt);
\filldraw [black] (3.3,3.2) circle (0.9pt);
\filldraw [black] (9.2,2.2) circle (0.9pt);
\filldraw [black] (8.9,3.2) circle (0.9pt);

\filldraw [black] (7.2,2) circle (0.9pt);

\draw (5,2) node[anchor=south]{$u$};
\draw (7.2,2) node[anchor=south]{$v$};
\draw (6.1,4) node[anchor=south]{$w$};
\end{tikzpicture}
\caption{The graph $\mathbf{C}_{(1):(1):(1,1)}$}
\label{fig}

\end{figure}
\end{Example}
\subsection{Technical lemma}

The above example shows that, for compact graphs of type~3, the
top-Betti contributions may fail to attain the upper bound. Therefore, when computing the multigraded Betti numbers of the edge
rings of such graphs, we also need to compute their second-top
multigraded Betti numbers. By Lemma~\ref{canonical}, this is
equivalent to computing the multigraded Betti numbers of the canonical
module of the edge ring in homological degree~$1$. The purpose of this subsection is to derive
an explicit formula for these Betti numbers.

Let $T$ be the subalgebra of the polynomial ring $\mathbb{K}[Y_1,\ldots,Y_m]$ generated by monomials $y_1,\ldots,y_n$. In \cite{V}, $T$ is called {\it a monomial algebra}. The monomials of $T$ form a monoid under multiplication, denoted by $\mathcal{U}$, and $T$ is isomorphic to the positive affine monoid ring $\mathbb{K}[\mathcal{U}]$. Let $R=\mathbb{K}[X_1,\ldots,X_n]$ be the polynomial ring, endowed with a $\mathbb{Z}^m$-grading by assigning the multidegree of $y_i$ to $X_i$.  Consider an ideal $I$ of $T$ generated by certain monomials from $\mathcal{U}$. By the existence of a surjective homomorphism from $R$ to $T$ mapping $X_i$ to $y_i$, it follows that $T$, $I$ and $T/I$ are all finitely generated $\mathbb{Z}^m$-graded $R$-modules.

A formula for computing the multigraded Betti numbers of $T/I$ was established in \cite[Proposition 1.1]{BH2} using relative simplicial homology. In this paper, however, we need to calculate the multigraded Betti numbers of $I$ as an $R$-module. For this, we associate   every multidegree $h\in \mathbb{Z}_{\geq0}^m$ with the following simplicial complex:
 $$\Gamma_h(I):=\{F\subseteq [n]\mid Y^h=y^F\alpha \mbox{ for some monomial } \alpha\in I\}.$$
Here,  $Y^h$ is  defined as $Y_1^{h_1}\cdots Y_m^{h_m}$ for $h=(h_1,\ldots,h_m)\in \mathbb{Z}_{\geq0}^m$, and $y^F$ represents the product of all monomials $y_i$ with $i\in F$.

\begin{Lemma}\label{formula}  Let $I$  be a monomial ideal of $T$, and let $h\in \mathbb{Z}_{\geq0}^m$ be a multidegree. Then, for all $i\geq 0$, one has  $$\beta_{i,h}^R(I)=\dim_\mathbb{K} \widetilde{H}_{i-1}(\Gamma_h(I), \mathbb{K}).$$
\end{Lemma}

\begin{proof} Let $K_{\bullet}$ denote the Koszul complex of $X_1,\ldots,X_n$. Namely, $K_{\bullet}$ is the following complex  $$0\rightarrow K_n\rightarrow\cdots\rightarrow K_{i}\rightarrow \cdots \rightarrow K_0\rightarrow 0. $$
Here, $K_i$ is a multigraded free $R$-module with basis $\{e_F\mid |F|=i,\ F\subseteq[n]\}$.  Moreover, $$\deg(e_F)=\sum_{i\in F}\deg(X_i)=\sum_{i\in F}\deg(y_i)\in \mathbb{Z}_{\geq0}^{m}.$$ In particular, $K_0=Re_{\emptyset}=R$.
Since $K_{\bullet}$ is a multigraded minimal free resolution of $\mathbb{K}$, we have $$\beta_{i,h}^R(I)=\dim_{\mathbb{K}}\mathrm{Tor}_{i}^R(\mathbb{K}, I)_{h}=\dim_{\mathbb{K}}H_i(K_{\bullet}\otimes_R I)_h.$$
On the other hand, $$(K_i\otimes I)_h=\bigoplus_{|F|=i}(Re_F\otimes_RI)_h=\bigoplus_{|F|=i}(I(\deg(e_F)))_h.$$
Note that $$(I(\deg(e_F)))_h=I_{h-\deg(e_F)}=\left\{
                                             \begin{array}{ll}
                                               \mathbb{K}, & \hbox{$F\in \Gamma_h(I)$;} \\
                                               0, & \hbox{otherwise.}
                                             \end{array}
                                           \right.$$
It follows that there is an isomorphism between the complex
$(K_{\bullet}\otimes_R I)_h$ and the complex
$\widetilde{C}_{\bullet}(\Gamma_h(I), \mathbb{K})[-1]$,
yielding the desired formula.
\end{proof}

We may also regard this lemma as a generalization of \cite[Theorem 1.34]{MS}, which addresses the same question on   monomial ideals in a polynomial ring.

\subsection{Type 3}
To apply Lemma~\ref{formula}, we first make an observation concerning the simplicial complexes $\Gamma_h(I)$.
Throughout this discussion, we use the notation of Lemma~\ref{formula}.
Let $\Omega$ be a set of monomials in $T$, and assume that the ideal $I$
is generated by the monomials in $\Omega$. For each $\alpha\in\Omega$, define
\[
\Gamma_h^{\alpha}(I)
=
\left\langle
\{i\in[n]\mid a_i\neq 0\}
\ \middle|\
Y^h=\alpha y_1^{a_1}\cdots y_n^{a_n},
\ a_1,\ldots,a_n\in\mathbb{Z}_{\geq 0}
\right\rangle,
\]
where $\langle\mathcal{F}\rangle$ denotes the simplicial complex generated
by the subsets in $\mathcal{F}$. If no such factorization exists, we set
$\Gamma_h^{\alpha}(I)=\emptyset$. Thus, the facets of
$\Gamma_h^{\alpha}(I)$ are precisely the inclusion-maximal members among
the supports arising from such factorizations.

Since $I$ is generated by the monomials in $\Omega$, it follows that
\[
\Gamma_h(I)
=
\bigcup_{\alpha\in\Omega}\Gamma_h^{\alpha}(I).
\]

For compact graphs of type~3, we will repeatedly use the following
observation. Set
\[
y_1=uv,\qquad y_2=vw,\qquad y_3=wu.
\]
Any factorization into edge monomials of a monomial involving only
$u,v,w$ can involve only $y_1,y_2,y_3$. Moreover, such a factorization,
if it exists, is unique. Indeed, if
\[
y_1^{a_1}y_2^{a_2}y_3^{a_3}
=
y_1^{b_1}y_2^{b_2}y_3^{b_3},
\]
then comparison of the exponents of $u,v,w$ gives
\[
a_1+a_3=b_1+b_3,\qquad
a_1+a_2=b_1+b_2,\qquad
a_2+a_3=b_2+b_3,
\]
and hence $a_i=b_i$ for $i=1,2,3$.

\begin{Example}\label{example} \em Let $\mathbb{K}[G]$ be the edge ring of $G=\mathbf{C}_{(1):(1):(1,1)}$. Note that all the vertices of $G$ are labeled, as shown in Figure~\ref{fig}.  Thus, $\Theta_G=uvwx_1x_2w_1w_2u_1u_2v_1v_2$. In the following, we denote $\Theta_G$ by $\Theta$ for short. Moreover, in view of Lemma \ref{3.4}, the canonical module
	$\omega_{\mathbb{K}[G]}$, as an ideal of $\mathbb{K}[G]$, is generated by
	monomials $$\Theta w, \Theta w^3, \Theta u, \Theta v.$$
	The four monomials are denoted by $\mu_i$ for $i=1,\ldots,4$.
Put $\alpha_1=\Theta uvw$ and $\alpha_2=\Theta uvw^3$.

Since $\Theta uvw=\mu_1 y_1$, it follows that $\Gamma_{\alpha_1}^{\mu_1}(\omega_{\mathbb{K}[G]})=\langle\{1\}\rangle$;

Since $\mu_2$ does not divide $\Theta uvw$, it follows that $\Gamma_{\alpha_1}^{\mu_2}(\omega_{\mathbb{K}[G]})=\emptyset$;

Since $\Theta uvw=\mu_3 y_2$, it follows that $\Gamma_{\alpha_1}^{\mu_3}(\omega_{\mathbb{K}[G]})=\langle\{2\}\rangle$;

Since $\Theta uvw=\mu_4 y_3$, it follows that $\Gamma_{\alpha_1}^{\mu_4}(\omega_{\mathbb{K}[G]})=\langle\{3\}\rangle$.

Thus, $\Gamma_{\alpha_1}(\omega_{\mathbb{K}[G]})=\langle\{1\}, \{2\}, \{3\}\rangle$, and so $\dim_{\mathbb{K}}\widetilde{H}_0(\Gamma_{\alpha_1}(\omega_{\mathbb{K}[G]}),\mathbb{K})=2$.

Similarly, we have $\Gamma_{\alpha_2}^{\mu_1}(\omega_{\mathbb{K}[G]})=\langle\{2,3\}\rangle$, $\Gamma_{\alpha_2}^{\mu_2}(\omega_{\mathbb{K}[G]})=\langle\{1\}\rangle$ and $\Gamma_{\alpha_2}^{\mu_3}(\omega_{\mathbb{K}[G]})=\Gamma_{\alpha_2}^{\mu_4}=\emptyset$. Hence, $\Gamma_{\alpha_2}(\omega_{\mathbb{K}[G]})=\langle\{2,3\},\{1\}\rangle$ and it follows that $\dim_{\mathbb{K}}\widetilde{H}_0(\Gamma_{\alpha_2}(\omega_{\mathbb{K}[G]}),\mathbb{K})=1$.

Note that $D_G=\Theta^2u^2v^2w^4$. Since $D_G/\alpha_1=\Theta uvw^3$ and $D_G/\alpha_2=\Theta uvw$,  it follows from Lemma \ref{formula} as well as Lemma \ref{canonical} that  $$\beta_{3, \Theta uvw^3}(\mathbb{K}[G])=\beta_{1,\alpha_1}(\omega_{\mathbb{K}[G]})=2$$ and $$\beta_{3, \Theta uvw}(\mathbb{K}[G])=\beta_{1,\alpha_2}(\omega_{\mathbb{K}[G]})=1.$$
Since $\operatorname{pdim}(\mathbb{K}[G])=4$ and these multidegrees
have full support, $G$ is a proper second-Betti graph.
There are precisely four induced subgraphs of  $G$ that consist of four minimal cycles of $G$. These four induced subgraphs are induced on $V(G)\setminus \{x_1,x_2\}$ and  $V(G)\setminus \{w_1,w_2\}$, $V(G)\setminus \{u_1,u_2\}$ and $V(G)\setminus \{v_1,v_2\}$, which we denote by $H_1,H_2,H_3,$ and $H_4$,
respectively. According to Lemmas~\ref{5.3} and~\ref{5.5}, we have $$\mathcal{N}_{H_1}=\Theta_{H_1}\{uv^2w^2, u^2vw^2, u^2v^2w\}, \quad \mathcal{N}_{H_2}=\Theta_{H_2}\{uv^2w^2, u^2vw^2, u^2v^2w \}$$ and $$\mathcal{N}_{H_3}=\Theta_{H_3}\{v^2w, v^2w^3,  vw^4\}, \quad \mathcal{N}_{H_4}=\Theta_{H_4}\{u^2w, u^2w^3,  uw^4\}.$$
Here, for a set $A$ of monomials and  for a monomial $\alpha$, the notation $\alpha A$ denotes the set $\{\alpha m\mid m\in A\}$.

Combining the above facts, the sum of the Betti numbers $\beta_{3,\alpha}(\mathbb{K}[G])$ with $\mbox{supp}(\alpha)=V(G)$ is  greater than or equal to 3, and the sum of the Betti numbers $\beta_{3,\alpha}(\mathbb{K}[G])$ with $\mbox{supp}(\alpha)=V(H_i)$ is also greater than or equal to 3 for $i=1,2,3,4$. Therefore, it follows that $$15\leq \beta_{3}(\mathbb{K}[G])\leq \beta_{3}(\mathbb{K}[E(G)]/\mathrm{in}_{<}(I_G))=(2+1)\binom{5}{4}=15.$$
Hence, $\beta_{3}(\mathbb{K}[G])= \beta_{3}(\mathbb{K}[E(G)]/\mathrm{in}_{<}(I_G))$, and $\beta_{3,\alpha}(\mathbb{K}[G])$ for $\alpha\in \mathbb{Z}^{V(G)}$ can also be determined.

\end{Example}

The approach used in Example~\ref{example} extends to the general case.

\begin{Definition}\label{6.6}{\em
Let $H$ be the graph $\mathbf{C}_{\underline{p}:\underline{q}:\underline{r}}$
with its vertices labeled as in Section~5.3. We define
\[
\begin{aligned}
\mathcal{M}_H^1
  &:= \Theta_H \Big(
       \{u^{2\ell-1}v^{2n-1}w^{2k-1} \mid 1\le \ell\le m-1\}
       \sqcup
       \{u^{2m-1}v^{2\ell-1}w^{2k-1} \mid 1\le \ell\le n-1\} \\
  &\hphantom{:= \Theta_H \Big(}\;\sqcup\;
       \{u^{2m-1}v^{2n-1}w^{2\ell-1} \mid 1\le \ell\le k-1\}
     \Big),\\
\mathcal{M}_H^2
  &:= \{\Theta_H u^{2m-1}v^{2n-1}w^{2k-1}\}.
\end{aligned}
\]
}
\end{Definition}

\begin{Lemma} \label{p-1} Let $H=\mathbf{C}_{\underline{p}:\underline{q}:\underline{r}}$ be a compact graph of type 3, and let $p=m+n+k$ be the projective dimension of $\mathbb{K}[H]$. For every $\alpha\in \mathcal{M}^1_H\sqcup \mathcal{M}^2_H$, one has $$\beta_{p-1,\alpha}(\mathbb{K}[H])=\left\{
                                                                                    \begin{array}{ll}
                                                                                      1, & \hbox{$\alpha\in \mathcal{M}^1_H$;} \\
                                                                                      2, & \hbox{$\alpha\in \mathcal{M}^2_H$.}
                                                                                    \end{array}
                                                                                  \right.
$$
In particular, we have $\sum\limits_{\alpha}\beta_{p-1,\alpha}(\mathbb{K}[H])\geq p-1$, where the summation is over all $\alpha\in \mathrm{Mon}(\mathbb{K}[H])$ with $\mathrm{supp}(\alpha)=V(H)$.
\end{Lemma}

\begin{proof}
	Let $\omega$ denote the canonical module $\omega_{\mathbb{K}[H]}$.
	Note that $\omega$ is a monomial ideal of $\mathbb{K}[H]$, whose
	generators are given in Lemma~\ref{3.4}.
Put $\alpha=\Theta uvw$, where $\Theta$ is used as a shorthand
for
$\Theta_H$. Then $$\Gamma_{\alpha}^{\Theta u}(\omega)=\langle\{2\}\rangle, \Gamma_{\alpha}^{\Theta v}(\omega)=\langle\{3\}\rangle, \Gamma_{\alpha}^{\Theta w}(\omega)=\langle\{1\}\rangle.$$ For other generators $\beta$ of $\omega$, one easily sees that $\Gamma_{\alpha}^{\beta}(\omega)=\emptyset$. This implies $\Gamma_{\alpha}(\omega)=\langle\{2\}, \{3\}, \{1\} \rangle$ and so $\beta_{1,\alpha}(\omega)=2$ by Lemma~\ref{formula}. In the same vein, we can prove $$\beta_{1,\Theta u^{2\ell-1}vw}(\omega)=1 \mbox{ for }  \ell=2,3,\ldots,m,$$
$$\beta_{1,\Theta uv^{2\ell-1}w}(\omega)=1 \mbox{ for }  \ell=2,3,\ldots,n,$$ and
$$\beta_{1,\Theta uvw^{2\ell-1} }(\omega)=1 \mbox{ for }  \ell=2,3,\ldots,k.$$

Note that $D_H=\Theta^2u^{2m}v^{2n}w^{2k}$. The first assertion follows from Lemma~\ref{canonical}.  For the second assertion, it follows from the facts that $|\mathcal{M}^1_H|=m+n+k-3$, $|\mathcal{M}^2_H|=1$,  and $\mathrm{supp}(\alpha)=V(H)$ for each $\alpha\in  \mathcal{M}^1_H\sqcup \mathcal{M}^2_H$.
\end{proof}

\begin{Proposition}\label{second} Let $G= \mathbf{C}_{\underline{p}:\underline{q}:\underline{r}}$ denote a compact graph of type 3. Then $$\beta_{i}(\mathbb{K}[G])=\beta_{i}(\mathbb{K}[E(G)]/\mathrm{in}_{<}(I_G))$$ for all $i\geq 0$.

\end{Proposition}
\begin{proof}
	We classify the minimal cycles of $G$ into four distinct classes. The cycle encompassing the three big vertices is designated as the big cycle. A minimal cycle that is distinct from the big cycle but includes the vertex $u$ (or $v$, or $w$) is referred to as a cycle containing $u$ (or $v$, or $w$ respectively). Thus, besides the big cycle, there are $m$ minimal cycles containing $u$, $n$ minimal cycles containing $v$, and $k$ minimal cycles containing $w$.
	
For $i=0$, both Betti numbers are equal to $1$. By Lemma~\ref{compactupper}, both sides are zero for $i>m+n+k$. Fix $1\leq i\leq m+n+k$.
	
	There are
	\[
	\binom{m+n+k+1}{i+1}
	\]
	choices of $i+1$ minimal cycles of $G$. Fix such a choice $C_1,\ldots,C_{i+1}$, and let $H$ be the induced subgraph of $G$ on
	\[
	\bigcup_{j=1}^{i+1}V(C_j).
	\]
	Every minimal cycle of $H$ is a minimal cycle of $G$. Moreover, for each non-big minimal cycle, every vertex other than its big
	vertex lies on no other minimal cycle. Hence, no unchosen non-big minimal
	cycle can be contained in $H$. Therefore, the only possible additional
	minimal cycle of $H$ is the big cycle.
	
	Suppose that the big cycle is not among $C_1,\ldots,C_{i+1}$. Since $H$ is induced, it contains the big cycle if and only if $u,v,w\in V(H)$, equivalently, if and only if the chosen cycles include at least one minimal cycle containing each of $u,v,w$. Thus, the choices for which an additional big cycle occurs are precisely those consisting of $a$ minimal cycles containing $u$, $b$ minimal cycles containing $v$, and $c$ minimal cycles containing $w$, where
	\[
	a+b+c=i+1,\qquad a,b,c\geq1.
	\]
	Their number is
	\[
	\sum\limits_{a+b+c=i+1}^{a,b,c\geq1}
	\binom{m}{a}\binom{n}{b}\binom{k}{c}.
	\]
	For each such choice, $H$ is a compact induced subgraph of type~3 containing exactly $i+2$ minimal cycles. For every remaining choice, $H$ is a compact induced subgraph containing exactly $i+1$ minimal cycles. In either case, the chosen cycles are uniquely determined by $H$, so distinct choices determine distinct induced subgraphs.
	
Consequently, we obtain
	\[
	\binom{m+n+k+1}{i+1}
	-
	\sum\limits_{a+b+c=i+1}^{a,b,c\geq1}
	\binom{m}{a}\binom{n}{b}\binom{k}{c}
	\]
	distinct compact induced subgraphs of $G$ containing exactly $i+1$ minimal cycles, and
	\[
	\sum\limits_{a+b+c=i+1}^{a,b,c\geq1}
	\binom{m}{a}\binom{n}{b}\binom{k}{c}
	\]
	distinct compact induced subgraphs of $G$ of type~3 containing exactly $i+2$ minimal cycles.
	
The multidegrees contributed by distinct induced subgraphs are
disjoint, since each such multidegree has support equal to the vertex
set of the corresponding induced subgraph.

Since every compact graph is top-Betti and
$|\mathcal{N}_H|=\mathfrak{t}(H)-1$, Lemma~\ref{start} gives
	\[
	\sum_{\alpha}\beta_{i,\alpha}(\mathbb{K}[G])
	\geq
	\left(
	\binom{m+n+k+1}{i+1}
	-
	\sum\limits_{a+b+c=i+1}^{a,b,c\geq1}
	\binom{m}{a}\binom{n}{b}\binom{k}{c}
	\right)i,
	\]
	where the summation is over all $\alpha\in\mathrm{Mon}(\mathbb{K}[G])$ such that the induced subgraph on $\mathrm{supp}(\alpha)$ is a compact graph containing exactly $i+1$ minimal cycles.
	
	On the other hand, Lemmas~\ref{p-1} and~\ref{start} give
	\[
	\sum_{\alpha}\beta_{i,\alpha}(\mathbb{K}[G])
	\geq
	\sum\limits_{a+b+c=i+1}^{a,b,c\geq1}
	\binom{m}{a}\binom{n}{b}\binom{k}{c}i,
	\]
	where the summation is over all $\alpha\in\mathrm{Mon}(\mathbb{K}[G])$ such that the induced subgraph on $\mathrm{supp}(\alpha)$ is a compact graph of type~3 containing exactly $i+2$ minimal cycles.
	
Since the two sums above are indexed by disjoint sets of monomials, it follows that
\[
\beta_i(\mathbb{K}[G])
\geq
\binom{m+n+k+1}{i+1}i.
\]
Now, the result follows from Lemma~\ref{compactupper}.\end{proof}
From the proof of Proposition~\ref{second}, we see that the inequality in Lemma~\ref{p-1} is actually an equality. Combining the preceding results, we immediately obtain the main result of this section.

\begin{Theorem} \label{main1} Let $G$ be a compact graph, and denote by $p = \mathfrak{t}(G) - 1$ the projective dimension of $\mathbb{K}[G]$. Consider an integer $i$ satisfying $1 \leq i \leq p$. Let $H_j$ for $j = 1, \ldots, s$ represent all the compact induced subgraphs of $G$ having exactly $i + 1$ minimal cycles, and let $L_j$ for $j = 1, \ldots, r$ be all the  compact induced subgraphs of $G$ belonging to type 3 with $i+ 2$ minimal cycles. Define the following sets:

$$
\mathcal{N}_i(G) = \bigsqcup_{j=1}^{s} \mathcal{N}_{H_j}, \quad \mathcal{M}_i^1(G) = \bigsqcup_{j=1}^{r} \mathcal{M}_{L_j}^1, \quad \text{and} \quad \mathcal{M}_i^2(G) = \bigsqcup_{j=1}^{r} \mathcal{M}_{L_j}^2.
$$
Then, the $i$-th multigraded Betti numbers of $\mathbb{K}[G]$ can be expressed as follows:

$$
\beta_{i,\alpha}(\mathbb{K}[G]) = \begin{cases}
2, & \text{if } \alpha \in \mathcal{M}_i^2(G), \\
1, & \text{if } \alpha \in \mathcal{M}_i^1(G) \sqcup \mathcal{N}_i(G), \\
0, & \text{otherwise}.
\end{cases}
$$
Note that $r$ is always equal to 0 if $G$ is a compact graph of type 1 or 2.
\end{Theorem}

For a compact graph $G$, let
$p=\operatorname{pdim}(\mathbb{K}[G])$. For $0\leq i\leq p-1$
and $j\geq0$, set
\[
O_{i,j}
=
\{\alpha\in \mathcal{M}^1_{i+1}(G)\sqcup \mathcal{N}_{i+1}(G)
\mid |\deg(\alpha)|=2j\}
\]
and
\[
T_{i,j}
=
\{\alpha\in \mathcal{M}^2_{i+1}(G)
\mid |\deg(\alpha)|=2j\}.
\]

\begin{Corollary}\label{main1co}
	With the above notation, let $I_G$ denote the toric ideal of $G$. There exists a squarefree initial ideal $\mathrm{in}_{<}(I_G)$ of $I_G$ such that
	\[
	\beta_{i,j}(\mathrm{in}_{<}(I_G))
	=
	\beta_{i,j}(I_G)
	=
	|O_{i,j}|+2|T_{i,j}|
	\]
	for all $0\leq i\leq p-1$ and $j\geq0$.
\end{Corollary}

In the following example,  to illustrate Theorem~\ref{main1}, we present the multigraded
Betti numbers (in homological degree $2$)  of the edge ring associated with the simplest type-3 compact graph.
\begin{Example}\label{ex:compact-betti}\em
Let $G$ be the compact graph $\mathbf{C}_{(1):(1):(1)}$.  For simplicity, we label its vertices as in Subsection~5.3 with a minor modification: $u_1, u_2$ instead of $u_{1,1}, u_{1,2}$, $v_1, v_2$ instead of $v_{1,1}, v_{1,2}$, and $w_1, w_2$ instead of $w_{1,1}, w_{1,2}$, as shown in Figure~5.

\begin{figure}[ht]
\centering
\begin{tikzpicture}[
  scale=1.0,
  bigv/.style={circle, fill=black, inner sep=1.5pt},
  smallv/.style={circle, fill=black, inner sep=0.9pt}
]

\begin{scope}[shift={(0,0)}]
  \node[bigv] (U1) at (0,0) {};
  \node[bigv] (V1) at (2.4,0) {};
  \node[bigv] (W1) at (1.2,2.0) {};

  \node[draw=none, fill=none] at (0,-0.35) {$u$};
  \node[draw=none, fill=none] at (2.4,-0.35) {$v$};
  \node[draw=none, fill=none] at (1.55,2.0) {$w$};

  \draw (U1)--(V1)--(W1)--(U1);

  \node[smallv] (U1a) at (-0.9,0.8) {};
  \node[smallv] (U1b) at (-0.9,-0.8) {};
  \draw (U1)--(U1a)--(U1b)--(U1);
  \node[draw=none, fill=none] at (-1.3,0.9) {$u_1$};
  \node[draw=none, fill=none] at (-1.3,-0.9) {$u_2$};

  \node[smallv] (V1a) at (3.3,0.8) {};
  \node[smallv] (V1b) at (3.3,-0.8) {};
  \draw (V1)--(V1a)--(V1b)--(V1);
  \node[draw=none, fill=none] at (3.7,0.9) {$v_1$};
  \node[draw=none, fill=none] at (3.7,-0.9) {$v_2$};

  \node[smallv] (W1a) at (0.6,3.0) {};
  \node[smallv] (W1b) at (1.8,3.0) {};
  \draw (W1)--(W1a)--(W1b)--(W1);
  \node[draw=none, fill=none] at (0.3,3.2) {$w_1$};
  \node[draw=none, fill=none] at (2.1,3.2) {$w_2$};

  \node[draw=none, fill=none] at (1.2,-1.6) {$G=\mathbf{C}_{(1):(1):(1)}$};
\end{scope}

\begin{scope}[shift={(5.5,0)}]
  \node[bigv] (U2) at (0,0) {};
  \node[bigv] (V2) at (2.4,0) {};
  \node[bigv] (W2) at (1.2,2.0) {};

  \node[draw=none, fill=none] at (0,-0.35) {$u$};
  \node[draw=none, fill=none] at (2.4,-0.35) {$v$};
  \node[draw=none, fill=none] at (1.55,2.0) {$w$};

  \draw (U2)--(V2)--(W2)--(U2);

  \node[smallv] (U2a) at (-0.9,0.8) {};
  \node[smallv] (U2b) at (-0.9,-0.8) {};
  \draw[dashed] (U2)--(U2a)--(U2b)--(U2);
  \node[draw=none, fill=none] at (-1.3,0.9) {$u_1$};
  \node[draw=none, fill=none] at (-1.3,-0.9) {$u_2$};

  \node[smallv] (V2a) at (3.3,0.8) {};
  \node[smallv] (V2b) at (3.3,-0.8) {};
  \draw (V2)--(V2a)--(V2b)--(V2);
  \node[draw=none, fill=none] at (3.7,0.9) {$v_1$};
  \node[draw=none, fill=none] at (3.7,-0.9) {$v_2$};

  \node[smallv] (W2a) at (0.6,3.0) {};
  \node[smallv] (W2b) at (1.8,3.0) {};
  \draw (W2)--(W2a)--(W2b)--(W2);
  \node[draw=none, fill=none] at (0.3,3.2) {$w_1$};
  \node[draw=none, fill=none] at (2.1,3.2) {$w_2$};

  \node[draw=none, fill=none] at (1.2,-1.6) {$H_1$};
\end{scope}

\begin{scope}[shift={(0,-5.4)}]
  \node[bigv] (U3) at (0,0) {};
  \node[bigv] (V3) at (2.4,0) {};
  \node[bigv] (W3) at (1.2,2.0) {};

  \node[draw=none, fill=none] at (0,-0.35) {$u$};
  \node[draw=none, fill=none] at (2.4,-0.35) {$v$};
  \node[draw=none, fill=none] at (1.55,2.0) {$w$};

  \draw (U3)--(V3)--(W3)--(U3);

  \node[smallv] (U3a) at (-0.9,0.8) {};
  \node[smallv] (U3b) at (-0.9,-0.8) {};
  \draw (U3)--(U3a)--(U3b)--(U3);
  \node[draw=none, fill=none] at (-1.3,0.9) {$u_1$};
  \node[draw=none, fill=none] at (-1.3,-0.9) {$u_2$};

  \node[smallv] (V3a) at (3.3,0.8) {};
  \node[smallv] (V3b) at (3.3,-0.8) {};
  \draw[dashed] (V3)--(V3a)--(V3b)--(V3);
  \node[draw=none, fill=none] at (3.7,0.9) {$v_1$};
  \node[draw=none, fill=none] at (3.7,-0.9) {$v_2$};

  \node[smallv] (W3a) at (0.6,3.0) {};
  \node[smallv] (W3b) at (1.8,3.0) {};
  \draw (W3)--(W3a)--(W3b)--(W3);
  \node[draw=none, fill=none] at (0.3,3.2) {$w_1$};
  \node[draw=none, fill=none] at (2.1,3.2) {$w_2$};

  \node[draw=none, fill=none] at (1.2,-1.6) {$H_2$};
\end{scope}

\begin{scope}[shift={(5.5,-5.4)}]
  \node[bigv] (U4) at (0,0) {};
  \node[bigv] (V4) at (2.4,0) {};
  \node[bigv] (W4) at (1.2,2.0) {};

  \node[draw=none, fill=none] at (0,-0.35) {$u$};
  \node[draw=none, fill=none] at (2.4,-0.35) {$v$};
  \node[draw=none, fill=none] at (1.55,2.0) {$w$};

  \draw (U4)--(V4)--(W4)--(U4);

  \node[smallv] (U4a) at (-0.9,0.8) {};
  \node[smallv] (U4b) at (-0.9,-0.8) {};
  \draw (U4)--(U4a)--(U4b)--(U4);
  \node[draw=none, fill=none] at (-1.3,0.9) {$u_1$};
  \node[draw=none, fill=none] at (-1.3,-0.9) {$u_2$};

  \node[smallv] (V4a) at (3.3,0.8) {};
  \node[smallv] (V4b) at (3.3,-0.8) {};
  \draw (V4)--(V4a)--(V4b)--(V4);
  \node[draw=none, fill=none] at (3.7,0.9) {$v_1$};
  \node[draw=none, fill=none] at (3.7,-0.9) {$v_2$};

  \node[smallv] (W4a) at (0.6,3.0) {};
  \node[smallv] (W4b) at (1.8,3.0) {};
  \draw[dashed] (W4)--(W4a)--(W4b)--(W4);
  \node[draw=none, fill=none] at (0.3,3.2) {$w_1$};
  \node[draw=none, fill=none] at (2.1,3.2) {$w_2$};

  \node[draw=none, fill=none] at (1.2,-1.6) {$H_3$};
\end{scope}

\end{tikzpicture}
\caption{The induced subgraphs $G, H_1,H_2,H_3$ of $\mathbf{C}_{(1):(1):(1)}$.}
\end{figure}

By Theorem~\ref{main1}, it suffices to consider compact induced subgraphs of $G$
with $4$ and $3$ minimal cycles.  Clearly $G$ itself is the unique compact
induced subgraph with $4$ minimal cycles.  The compact induced subgraphs with
three minimal cycles are precisely the following three:
$H_1 := G\setminus\{u_1,u_2\},\quad
H_2 := G\setminus\{v_1,v_2\},\quad
H_3 := G\setminus\{w_1,w_2\}.$

By Definition~\ref{6.6} we have
\[
\mathcal{M}_2^1(G)=\mathcal{M}_G^1 = \emptyset,
\qquad
\mathcal{M}_2^2(G) = \mathcal{M}_G^2=\{\Theta_G\,uvw\}
                 = \{u_1u_2v_1v_2w_1w_2u^2v^2w^2\}.
\]
Moreover, by Lemma~\ref{5.3}(2) we obtain
\[
\begin{aligned}
\mathcal{N}_2(G)
  &= \mathcal{N}_{H_1}\sqcup\mathcal{N}_{H_2}\sqcup\mathcal{N}_{H_3} \\
  &= \{\Theta_{H_1}v^2w,\ \Theta_{H_1}vw^2,\
      \Theta_{H_2}u^2w,\ \Theta_{H_2}uw^2,\
      \Theta_{H_3}u^2v,\ \Theta_{H_3}uv^2\} \\
  &= \{v_1v_2w_1w_2uv^3w^2,\ v_1v_2w_1w_2uv^2w^3,\
      u_1u_2w_1w_2u^3vw^2,\\
  &\qquad
      u_1u_2w_1w_2u^2vw^3,\
      u_1u_2v_1v_2u^3v^2w,\
      u_1u_2v_1v_2u^2v^3w\}.
\end{aligned}
\]

Applying Theorem~\ref{main1} and Corollary~\ref{main1co} to $G=\mathbf{C}_{(1):(1):(1)}$
in homological degree $2$, we obtain
\[
\beta_{2,\alpha}(\mathbb{K}[G])=
\begin{cases}
  2, & \text{if }\alpha = u_1u_2v_1v_2w_1w_2u^2v^2w^2,\\[2pt]
  1, & \text{if }\alpha \in \mathcal{N}_2(G),\\[2pt]
  0, & \text{otherwise,}
\end{cases}
\qquad
\beta_{2,j}(\mathbb{K}[G])=
\begin{cases}
  6, & \text{if } j=5,\\[2pt]
  2, & \text{if } j=6,\\[2pt]
  0, & \text{otherwise.}
\end{cases}
\]
\end{Example}

{\bf \noindent Acknowledgment.} The authors would like to thank the anonymous referees for their careful reading
and helpful comments, which improved the presentation of this paper.
The authors are grateful to the computer algebra systems
CoCoA~\cite{C} and Macaulay2~\cite{GS} for providing a large number of examples
that helped us develop ideas and test our results.

\medskip

{\bf \noindent Declarations.}

\noindent\textbf{Funding.}
This work was supported by the National Natural Science Foundation of China (Grant No.~11971338).

\medskip

\noindent\textbf{Competing interests.}
The authors have no competing interests to declare that are relevant to the content of this article.

\medskip

\noindent\textbf{Data availability.}
No data are associated with this work.

\end{document}